\documentclass{amsart}
\usepackage[vcentermath,enableskew]{youngtab}
\usepackage{amsmath}
\usepackage{amsopn}
\usepackage{amssymb}
\usepackage{latexsym}
\usepackage{amsfonts}
\usepackage{amsthm}

\newcommand{\K}{{\mathcal K}}
\newcommand{\T}{{\rho}}
\newcommand{\hk}{{h}}

\newcommand{\der}{{\Delta\,}}
\newcommand{\ua}{{\dot{a}}}
\newcommand{\ub}{{\dot{b}}}
\newcommand{\da}{{\ddot{a}}}
\newcommand{\db}{{\ddot{b}}}

\newcommand{\then}{\Longrightarrow}

\newcommand{\sn}{{\it span}}

\newcommand{\al}{\alpha}
\newcommand{\G}{{X}}

\newtheorem{thm}{Theorem}[section]

\newtheorem{pro}[thm]{Proposition}
\newtheorem{lem}[thm]{Lemma}
\newtheorem{cla}[thm]{Claim}

\newtheorem{cor}[thm]{Corollary}

\theoremstyle{definition}
\newtheorem{obs}[thm]{Observation}
\newtheorem{fac}[thm]{Fact}
\newtheorem{rem}[thm]{Remark}
\newtheorem{exa}[thm]{Example}
\newtheorem{defn}[thm]{Definition}
\newtheorem{df}[thm]{Definition}

\newtheorem{conj}[thm]{Conjecture}

\newcommand{\een}{\end{enumerate}}
\newcommand{\blem}{\begin{lem}}
\newcommand{\elem}{\end{lem}}
\newcommand{\bcl}{\begin{cla}}
\newcommand{\ecl}{\end{cla}}
\newcommand{\ethm}{\end{thm}}
\newcommand{\bpr}{\begin{pro}}
\newcommand{\epr}{\end{pro}}
\newcommand{\bco}{\begin{cor}}
\newcommand{\eco}{\end{cor}}
\newcommand{\bcon}{\begin{conj}}
\newcommand{\econ}{\end{conj}}
\newcommand{\bde}{\begin{defn}}
\newcommand{\ede}{\end{defn}}
\newcommand{\bex}{\begin{exa}}
\newcommand{\eexa}{\end{exa}}
\newcommand{\bobs}{\begin{obs}}
\newcommand{\eobs}{\end{obs}}
\newcommand{\bexe}{\begin{exe}}
\newcommand{\eexe}{\end{exe}}

\newcommand{\la}{\lambda}

\newcommand{\mone}{{-1}}
\newcommand{\pmone}{{\pm1}}
\newcommand{\pmtwo}{{\pm2}}
\newcommand{\mtwo}{{-2}}
\newcommand{\mthree}{{-3}}
\newcommand{\pmthree}{{\pm3}}
\newcommand{\mfour}{{-4}}

\newcommand{\mfive}{{-5}}
\newcommand{\msix}{{-6}}
\newcommand{\mseven}{{-7}}
\newcommand{\meight}{{-8}}
\newcommand{\mnine}{{-9}}

\newcommand{\ma}{{-a}}
\newcommand{\mb}{{-b}}
\newcommand{\mc}{{-c}}

\newcommand{\plmb}{{\pm b}}

\begin{document}
\title{Abstract Young Pairs for Signed Permutation Groups }
\author{Yona Cherniavsky}
\address{Department of Mathematics and Statistics
 Bar-Ilan University
 Ramat-Gan, Israel 52900}
\email{cherniy@math.biu.ac.il , chrnvsk@gmail.com}
\thanks{Supported in part by a grant from 
the Israel Science Foundation this work is a part of the author's PhD research carried out at Bar Ilan University, Ramat Gan, Israel.}
\begin{abstract}
The notion of an Abstract Young (briefly: AY) representation is a
natural generalization of the classical Young orthogonal form. The
AY representations of the symmetric group are characterized in
\cite{U-II}. In this paper we present several types of minimal AY
representation of $D_n$ associated with standard D-Young tableaux
which are a natural generalization of usual standard Young
tableaux. We give an explicit combinatorial view (the
representation space is spanned by certain standard tableaux while
the action is a generalized Young orthogonal form) of
representations which are induced into $D_n$ from minimal AY
representations of one of the natural embeddings of $S_n$ into
$D_n$. Then we show that these induced representations are
isomorphic to the direct sum of two or three minimal AY
representations of $D_n$ also associated with standard D-Young
tableaux. It is done by constructing a continuous path between
representation matrices where one end of the path is the mentioned
direct sum; another end is the classical form of induced
representation. In the last  section we briefly explain how the
similar results may be obtained for the group $B_n$ instead of
$D_n$.
\end{abstract}
\maketitle

\section{Introduction}

 One of the important problems in combinatorial representation
theory is to find a unified combinatorial construction of Weyl
groups representations. A most important breakthrough in this area
wasthe introduction of Kazhdan-Lustig cell representations
\cite{KL}. Other results were achieved by Vershik \cite{V},
Vershik and Okounkov \cite{OV}, Cherednik \cite{Cher} and Ram
\cite{Ra}.

A unified axiomatic approach, in steps of the above works of
Vershik and Ram,
  to the representation theory of
Coxeter groups and their Hecke algebras was presented
in~\cite{U-I}. This was carried out by a natural assumption on the
representation matrices, avoiding a priori use of  external
concepts (such as Young tableaux).

\smallskip

Let $(W,S)$ be a Coxeter system, and let $\K$ be a finite subset
of $W$. Let $\mathbb F$ be a suitable field of characteristic zero
and let $\T$ be a representation of $W$ on the vector space
$V_{\K}:=\sn_{\mathbb F}\{C_w\,|\,w\in\K\}$, with basis vectors
indexed by elements of $\K$. Adin, Brenti and Roichman in
\cite{U-I} and \cite{U-II} study the sets $\K$ and representations
$\T$ which satisfy the following axiom:

\begin{itemize}
\item[$(A)$] {\it For any generator $s\in S$ and any element $w\in
\K$ there exist scalars $a_s(w),b_s(w)\in\mathbb F$ such that
$$
\T_s(C_w) = a_{s}(w)C_w + b_{s}(w)C_{ws}.
$$
If $w\in\K$ but $ws\not\in\K$ we assume $b_s(w)=0$.%
}
\end{itemize}

A pair $(\T,\K)$ satisfying Axiom $(A)$ is called an {\em abstract
Young (AY) pair;} $\T$ is an {\em AY representation,} and $\K$ is
an {\em AY cell.} If $\K\ne \varnothing$ and has no proper subset
$\varnothing \subset \K' \subset \K$ such that $V_{\K'}$ is
$\T$-invariant, then $(\T,\K)$ is called a {\em minimal AY pair.}
(This is much weaker than assuming $\T$ to be irreducible.)

 In~\cite{U-I} it was shown that  an AY representation of a
simply laced Coxeter group is determined by a linear functional on
the root space. Thus it may be obtained by restriction of Ram's
calibrated representations of affine Hecke algebras (see
\cite{Ra}) to the corresponding Weyl groups. In \cite{U-II} it is
shown that, furthermore, the values of the linear functional on
the ``boundary" of the AY cell determine the representation.

In \cite{U-II} this result is used to characterize AY cells in the
symmetric group. This characterization is then applied to show
that every irreducible representation of $S_n$ may be realized as
a minimal abstract Young representation. AY representations of
Weyl groups of type $B$ are not determined by a linear functional.
However, it is shown in \cite{U-II} that irreducible
representations of $B_n$, similarly to irreducible representations
of $S_n$, may be realized as minimal AY representations.

In this work we present several types of minimal AY representation
of $D_n$ which arise from D-Young tableaux introduced in Section
\ref{s.IdCT}. These D-Young tableaux are a special case of Ram's
negative rotationally symmetric tableaux \cite{Ra}. It is shown in
{\rm \cite[Theorem 3.9]{U-II}} that the representation of a
Coxeter group $W$ which is induced from a minimal AY
representation of its parabolic subgroup $P$ is a minimal AY
representation of $W$. In section \ref{indSntoDn} of this work we
give an explicit combinatorial view (the representation space is
spanned by certain standard tableaux while the action is a
generalized Young orthogonal form) of representations which are
induced into $D_n$ from minimal AY representations of one of the
natural embeddings of $S_n$ into $D_n$. In section \ref{decind} we
show that these induced representations are isomorphic to the
direct sum of two or three minimal AY representations of $D_n$
also associated with standard D-Young tableaux. It is done by
constructing a continuous path between representation matrices
where one end of the path is the mentioned direct sum; another end
is the classical form of induced representation (see {\rm
\cite[Remark 3.10]{U-II}}). In section~\ref{bn} we briefly explain
how the similar results may be obtained for the group $B_n$
instead of $D_n$.

\section{Preliminaries and notations}\label{prenot}

A {\em Coxeter system} is a pair $(W,S)$ consisting of a group $W$
and a set $S$ of generators for $W$, subject only to relations of
the form
$$
(st)^{m(s,t)}=1,
$$
where $m(s,s)=1$ and $m(s,t)=m(t,s)\ge 2$ for $s\ne t$ in $S$. In
case no relation occurs for a pair $(s,t)$, we make the convention
that $m(s,t)=\infty$. $W$ is called a {\em Coxeter group}. If
$m(s,t)\le 3$ for all $s\ne t$ then $(W,S)$ is called {\em simply
laced}.

Let $(W,S)$ be a Coxeter system, and let $\G=\G(W,S)$ be the
corresponding Cayley graph (with generators acting on the right):
its vertices are the elements of $W$, and $x, y\in W$ are
connected by an edge if and only if $x^{-1}y\in S$. $\G$ is a
connected undirected graph. A subset $\K \subseteq W$ is called
{\em convex} if, for any $x, y\in \K$, all the geodesics (paths of
shortest length) connecting $x$ to $y$ in $\G$ have all their
vertices in $\K$.

Let $P$ be a poset, and let $\G$ be its undirected Hasse diagram.
Thus $\G$ is an undirected graph which has $P$ as a vertex set,
with an edge $\{x,y\}$ whenever $x$ either covers or is covered by
$y$. A subset $\K$ of $P$ is called {\em convex} if, for any
$x,y\in\K$, all geodesics (shortest paths) connecting $x$ to $y$
in $\G$ have all their vertices in $\K$.

Below, $P$ will be a Coxeter group $W$ with the right weak Bruhat
order, namely the transitive closure of the relation
$$
w<ws \iff w\in W, s\in S \hbox{\ and\ } \ell(w)<\ell(ws).
$$
Clearly, $\K\subseteq W$ is convex in the right Cayley graph
$\G(W,S)$ if and only if it is convex in the right weak Bruhat
poset $P$.

\bigskip

Let $V$ be the root space of a Coxeter system $(W,S)$.
Thus $V$ is a vector space over $\mathbb R$ with a basis
$\{\al_s\,|\,s\in S\}$ indexed by the group generators. A
symmetric bilinear form $B$ is defined on $V$ by
$$
B(\al_s,\al_t) := -\cos{\pi\over m(s,t)}\qquad(\forall s,t\in S)
$$
(interpreted to be $-1$ in case $m(s,t)=\infty$). For each
generator $s\in S$ define a linear map $\sigma_s:V \to V$ by
$$
\sigma_s(v) := v - 2B(v,\al_s)\al_s\qquad(\forall v\in V).
$$
This yields a faithful $B$-preserving action $\sigma$ of $W$ on
$V$ (see, e,g.,~\cite{Hum}), and defines the corresponding {\em
root system}
$$
\Phi := \{\sigma_w(\al_s)\,|\,w\in W, s\in S\} \subseteq V.
$$

\subsection{The Coxeter groups of type B}\label{bdef}
The Coxeter group of type $B$, $B_n$, is the group of all signed
permutations. Let $S_{2n}$ be the group of all permutations of the
numbers $\pm1,\pm2,...,\pm n$. Then
$$
B_n=\{\pi\in
S_{2n}\,\,:\,\,\pi(-i)=-\pi(i)\,\,\textrm{for}\,\,i=1,2,...,n\}
$$
The group $B_n$ is generated by the Coxeter generators
$\{s_0,s_1,..,s_{n-1}\}$ defined by:
$$s_0=(-1,1)$$ and
$$s_i=(i,i+1)(-i,-i-1),\qquad 1 \leq i \leq n-1.$$

\subsection{The Coxeter groups of type D}\label{ddef}

The Coxeter group of type $D$, $D_n$, can be defined as the normal
subgroup of $B_n$ consisting of all signed permutations $\pi$
satisfying: $|\{i \in [n] \mid \pi(i)<0\}|$ is even. We embed
$D_n$ in $S_{2n}$ in the natural way. The group $D_n$ is generated
by the Coxeter generators $\{s_0,s_1,..,s_{n-1}\}$ defined by:
$$s_0=(1,-2)(2,-1)$$ and
$$s_i=(i,i+1)(-i,-i-1),\qquad 1 \leq i \leq n-1.$$
Notice that $s_0$ of $B_n$ differs from $s_0$ of $D_n$ while other
Coxeter generators are the same. An element of $D_n$ is called a
{\it reflection} if it is conjugate to a Coxeter generator. The
reflections of $D_n$ are $(i,j)=(i,j)(-i,-j)\in S_{2n}$. Notice
that in this notation $(i,j)=(-i,-j)$.

Let $V=\mathbb{R}^n$ be the root space of $D_n$ with
$\{e_1,\ldots,e_n\}$ as its standard basis. The simple roots are
$-e_1-e_2,e_1-e_2, e_2-e_3,\ldots,e_{n-1}-e_n$. Denote
$$
e_{-i}=-e_i\quad\textrm{for}\quad 1\leq i\leq n
$$
The positive root $\alpha_{ij}\in \mathbb R^n$ corresponding to
the reflection $(i,j)\in D_n$ is
$$
\alpha_{ij}= e_i-e_j\,\,\,\textrm{for}\quad i,j\in [\pm
n]\,\,,\,\,i\neq-j\,\,,\,\,i<j
$$
Notice that $\alpha_{ij}=\alpha_{-j,-i}=-\alpha_{ji}$. We have the
following identity:
$$
\alpha_{ir}+\alpha_{rj}=e_i-e_r+e_r-e_j=e_i-e_j=\alpha_{ij}
$$

\bde \label{dervect} For $v=(v_1,v_2,\ldots,v_n)\in\mathbb R^n$
denote $v_{-i}=-v_i$ and define the derived vector
\begin{align*}
\Delta
v&=(v_2-v_{-1},v_2-v_1,v_3-v_2,\ldots,v_n-v_{n-1})\\
&=(v_1+v_2,v_2-v_1,v_3-v_2,\dots,v_n-v_{n-1})\in\mathbb R^n\,\,.
\end{align*}
\ede

\subsection{The Young Orthogonal form}
Let $\la$ be a partition of $n$, and let $\T^\la$ be the
corresponding irreducible representation of $S_n$ (on the Specht
module $S^\la$). Let $\{v_Q\,|\,Q$ standard Young tableau of shape
$\la\}$ be the basis of $S^\la$ obtained from the basis of
polytabloids by the Gram-Schmidt process (see~\cite[\S 25]{Ja}).
Let $Q$ be a standard Young tableau. If $k\in \{1,\dots,n\}$ is in
box $(i,j)$ of $Q$, then the {\em content} of $k$ in $Q$ is
$$
c_k:= j - i.
$$
The $k$th {\em hook-distance} is defined as
$$
\hk(k):=c_{k+1}-c_k\qquad(1\le k\le n-1).
$$
Finally, let $Q^{s_i}$ be the tableau obtained from $Q$ by
interchanging $i$ and $i+1$.

\begin{thm}\label{t.yof}{\bf (Classical Young Orthogonal Form)}
{\rm\cite[\S 25.4]{Ja}} In the above notations,
$$
\T^\la(s_i)(v_Q)= \frac{1}{\hk(i)}\,v_Q +
\sqrt{1-\frac{1}{\hk(i)^2}}\,v_{Q^{s_i}}.
$$
\end{thm}

This setting is generalized naturally to skew shapes $\la/\mu$.
Standard Young tableaux of shape $\la/\mu$, hook distances on
these tableaux and skew Specht modules $S^{\la/\mu}$ are defined
analogously. In particular,

\begin{thm}\label{t.yof-skew}
{\bf (Classical Young Orthogonal Form for Skew Specht Modules)}
Let $\{v_Q\,|\,Q$ standard Young tableau of shape $\la/\mu\}$ be
the basis of the skew Specht module $S^{\la/\mu}$ obtained by the
Gram-Schmidt process from the polytabloid basis. Then
\begin{equation}\label{e.yofi}
\T^{\la/\mu}(s_i)(v_Q)= \frac{1}{\hk(i)}\,v_Q +
\sqrt{1-\frac{1}{\hk(i)^2}}\,v_{Q^{s_i}}.
\end{equation}
\end{thm}

\subsection{Abstract Young Cells and
Representations}\label{s.survey}

This section surveys
results from~\cite{U-I} which will be used in this paper. Recall
the definition of AY cells and representations from the
introduction.

\begin{obs}\label{t.main_unitcell}{\rm\cite[Observation 3.3]{U-I}}
Every nonempty AY cell is a left translate of an AY cell
containing the identity element of $W$.
\end{obs}

\begin{pro}\label{t.main.convex}{\rm\cite[Corollary 4.4]{U-I}}
Every minimal AY cell is convex (in the right Cayley graph $\G(W,S)$
or, equivalently, under right weak Bruhat order). In particular,
$b_s(w)\ne 0$ whenever $s\in S$ and $w, ws\in\K$.
\end{pro}

Important examples of AY cells are (standard) descent classes.

\bigskip

Surprisingly, Axiom $(A)$ leads to very concrete matrices, whose
entries are essentially inverse linear.

In~\cite{U-I}
it is shown that, under mild conditions, Axiom $(A)$ is equivalent
to the following more specific version. Here $T$ is the set of all
reflections in $W$.

\begin{itemize}
\item[$(B)$]
{\it For any reflection $t\in T$ there exist scalars $\ua_t,
\ub_t, \da_t, \db_t \in\mathbb F$ such that, for all $s\in S$ and
$w\in \K$:

$$
\T_s(C_w) =\begin{cases} \ua_{wsw^{-1}}C_w+\ub_{wsw^{-1}}C_{ws},
&if \ell(w)<\ell(ws);\cr \da_{wsw^{-1}}C_w+\db_{wsw^{-1}}C_{ws},
&if \ell(w)>\ell(ws).\cr\end{cases}
$$
If $w\in\K$ and $ws\not\in\K$ we assume that $\ub_{wsw^{-1}}=0$
(if $\ell(w)<\ell(ws)$) or
$\db_{wsw^{-1}}=0$ (if $\ell(w)>\ell(ws)$).%
}
\end{itemize}

\begin{thm}\label{t.main.a1}{\rm\cite[Theorem 5.2]{U-I}}
Let $(\T,\K)$ be a minimal AY pair for the Iwahori-Hecke algebra
of $(W,S)$.
If 
$a_s(w)=a_{s'}(w')\then b_s(w)=b_{s'}(w')$ $(\forall s,s'\in S,
w,w'\in\K)$, then $\rho$ satisfies Axiom $(B)$.
\end{thm}

\begin{thm}\label{t.main.c1}{\rm\cite[Theorem 11.1]{U-I}}
The coefficients $\ua_t$ $(t\in T)$ determine all the character
values of $\T$.

\end{thm}

The assumption regarding the coefficients $b_s(w)$ in
Theorem~\ref{t.main.a1} is merely a normalization condition.
Theorem~\ref{t.main.a1} shows that the coefficients $a_s(w)$ and
$b_s(w)$ in Axiom $(A)$ depend only on the reflection $wsw^{-1}\in
T$ and on the relation between $w$ and $ws$ in the right weak
Bruhat order. It turns out that for simply laced Coxeter groups
the coefficients $\ua_t$ are given by a linear functional.


\begin{df}
For a convex subset $\K\subseteq W$ define:
\begin{eqnarray*}
T_{\K}          &:=& \{wsw^{-1}\,|\,s\in S,\,w\in \K,\,ws\in \K\},\\
T_{\partial \K} &:=& \{wsw^{-1}\,|\,s\in S,\,w\in \K,\,ws\not\in \K\}.
\end{eqnarray*}
\end{df}

\begin{df}\label{d.cg} ($\K$-genericity)\\
Let $\K$ be a convex subset of $W$ containing the identity
element. A vector $f$ in the root space $V$ is {\rm $\K$-generic}
if:
\begin{itemize}
\item[{\bf (i)}] For all $t\in T_{\K}$,
$$
\langle f,\al_t\rangle \not\in \{0, 1, -1\}.
$$
\item[{\bf (ii)}] For all $t\in T_{\partial \K}$,
$$
\langle f,\al_t\rangle \in \{1, -1\}.
$$
\item[{\bf (iii)}] If $w\in \K$, $s,t\in S$, $m(s,t)=3$ and $ws,
wt\not\in \K$, then
$$
\langle f,\al_{wsw^{-1}}\rangle = \langle
f,\al_{wtw^{-1}}\rangle\;(= \pm1).
$$
\end{itemize}
\end{df}

By Observation~\ref{t.main_unitcell}, we may assume that
$id\in\K$.

\begin{thm}\label{t.main11}{\rm\cite[Theorem 7.4]{U-I}}
Let $(W,S)$ be an irreducible simply laced Coxeter system, and let
$\K$ be a convex subset of $W$ containing the identity element. If
$f\in V^*$ is $\K$-generic, then
$$
\ua_t := {1\over \langle f,\al_t\rangle} \qquad(\forall t\in
T_{\K} \cup T_{\partial\K}),
$$
together with $\da_t$, $\ub_t$ and $\db_t$ satisfying
\begin{eqnarray*}
\ua_t + \da_t &=& 0\\
\ub_t \cdot \db_t &=& (1 - \ua_t)(1 - \da_t)
\end{eqnarray*}
define a representation $\T$ such that $(\T,\K)$ is a minimal AY
pair satisfying Axiom $(B)$.
\end{thm}

\begin{rem}\label{normalization} Various normalizations for $\ub_t$ and $\db_t$ are possible:
symmetric ($\db_t = \ub_t$), seminormal ($\db_t = 1$), row
stochastic ($\ua_t + \ub_t = \da_t + \db_t = 1$), etc.\
\cite[Subsection 5.2]{U-I}. By Theorem~\ref{t.main.c1}, all these
normalizations give isomorphic representations.\end{rem} The
following theorem is complementary.
\begin{thm}\label{t.main12}{\rm\cite[Theorem 7.5]{U-I}}
Let $(W,S)$ be an irreducible simply laced Coxeter system and let
$\K$ be a subset of $W$ containing the identity element. If
$(\T,\K)$ is a minimal AY pair satisfying Axiom $(B)$ and
$\ua_t\ne 0$ ($\forall t\in T_{\K}$), then there exists a
$\K$-generic $f\in V$ such that
$$
\ua_t = \frac{1}{\langle f,\al_t\rangle}\qquad (\forall\ t\in
T_{\K} \cup T_{\partial\K}).
$$
\end{thm}

\subsection{Boundary Conditions}\label{s.bound}

In this section it is shown that the action of the group $W$ on
the boundary of a minimal AY cell determines the representation up
to isomorphism.

\begin{df}\label{d.CF}
Let $f\in V$ be an arbitrary vector on the root space $V$ of $W$.
\begin{itemize}
\item[(1)] Define
$$
\K^f:=\{w \in W\mid \forall t \in A_f,  \ell(tw)>\ell(w)\}$$ where
$$A_f=\{t \in T \mid \langle f,\alpha_t \rangle \in \{\pm 1\}\}.$$
\item[(2)] If $f$ is $\K^f$-generic (as in Definition~\ref{d.cg}),
then the corresponding AY representation of $W$ (as in
Theorem~\ref{t.main11}), with the symmetric normalization $\db_t
=\ub_t\,(\forall t\in T_{\K^f})$, will be denoted $\T^f$.
\end{itemize}
\end{df}

\section{D-Young Tableaux and Minimal Cells in $D_n$}\label{s.cells_Sn}

In this section we show that standard Young tableaux of skew shape
lead to minimal AY cells.

\subsection{Cells and Skew Shapes}\label{s.IdCT}

In this subsection we study minimal AY cells $\K\subseteq D_n$. By
Observation~\ref{t.main_unitcell}, every minimal AY cell is a
translate of a minimal AY cell containing the identity element;
thus we may assume that $id\in\K$.

\medskip

We identify the root space of $D_n$ with $V=\mathbb R^{n}$ and for
a vector $v=(v_1,\dots,v_n)\in \mathbb R^n$ and recall the
notation from Definition \ref{dervect}
$$
\der v=(v_1+v_2,v_2-v_1,\dots,v_n-v_{n-1})\in \mathbb R^{n}.
$$

For a (skew) tableau $T$ denote $$c_k:=j-i$$ where $k$ is the
entry in row $i$ and column $j$ of $T$. Call
$cont(T):=(c_1,\dots,c_n)$ the {\em content vector} of $T$, and
call $\der cont(T)$ the {\em derived content vector} of $T$. Below
we sometimes shall denote for brevity $cont(T)$ as $c(T)$.

\bde\label{dyoungtab}
 Let $\lambda$ be a diagram of a skew shape. Define a $D$-Young
tableau of shape $\lambda$ to be a filling of $\lambda$ by the
$2n$ numbers $\pm 1,\pm 2,\ldots,\pm n$ in such a way that
$c_{-i}=-c_i$ for $1 \leq i \leq n$. A $D$-Young tableau is called
{\it standard} if the numbers are increasing in rows and in
columns. If $c_i=0$, then we allow the numbers $\pm i$ to occupy
the same box. \ede
\begin{rem}
Our standard D-Young tableau doesn't change if we multiply all its
entries by $-1$ and rotate it by $180^\circ$. Indeed standard
D-Young tableaux (except of the case when $\pm i$ occupy the same
box) are ``negative rotationally symmetric standard tableaux"
considered by Ram in~\cite{Ra} just with additional requirement
that the all assigned boxes form a skew shape.
\end{rem}
\begin{rem}\label{intcont}
\begin{itemize}
\item[(1)] The contents $c_i$ are \textbf{integer} numbers
...,$-2$, $-1$,$0$,$1$,$2$,... We emphasize it here because in
section~\ref{bn} we shall consider noninteger contents.
\item[(2)]
 The need of considering the negative entries and the negative contents, as well as the need to allow
 two numbers in the same box will be discussed
  in subsection~\ref{necessneg}.
\end{itemize}
\end{rem}
\Yautoscale0 \Yboxdim16pt
 Here are three simple examples:
$$\young(:\mnine\meight\mseven\msix\mfive,:\mfour\mthree\mtwo\mone,:1234,56789)
 \qquad  \qquad\young(:\mfour1,\mthree\pmtwo3,\mone4)  \qquad\qquad
 \young(::\mthree\mtwo14,\mfour\mone23)\quad.$$

\medskip

Recall the notations $\K^f$ and $\T^f$ from Definition~\ref{d.CF}.
Also recall that for each $v \in \mathbb R^n$, $\langle
v,\alpha_{ij}\rangle=v_i-v_j$.

The following theorem generalizes the part of sufficiency of
Theorem 5.1 from \cite{U-II}.

\begin{thm}\label{tablestoreps}

Let $T$ be a standard $D$-Young tableau with $2n$ boxes, let
$c=(c_1,...,c_n)\in\mathbb Z^n$ be its content vector. Consider
the sets $A_c=\{t\in T \mid \langle c,\alpha_t \rangle \in
\{1,-1\}\}$ and $\K^c=\{w \in D_n \mid \ell(tw)>\ell(w),\quad
\forall t \in A_c\}$. Then $c$ is $K^c$-generic and therefore
gives rise to a minimal AY pair $(\rho^c,\K^c)$.
\end{thm}
\noindent See~{\rm\cite[Theorem~2.3.4]{Cr}} for the proof of
Theorem~\ref{tablestoreps}.
\subsection{DAY Cells: definition and structure.}\label{s.CET}

\bde\label{DAY} Let $T$ be a standard $D$-Young tableau with $2n$
boxes. It follows from Theorem~\ref{tablestoreps}, that $T$ gives
rise to a minimal $AY$ pair $(\K^c,\rho^c)$ where $c=cont(T)$. We
call such a cell $K^c$ a DAY cell and such a representation
$\rho^c$ a DAY representation. \ede The following theorem
describes a DAY cell by a certain set of standard D-Young
tableaux. Its statement and proof are similar to the statement and
the proof of Theorem~5.5 from~\cite{U-II}. The proof of
Theorem~\ref{t.cell} may be found
in~{\rm\cite[Theorem~2.3.12]{Cr}}.
\begin{thm}\label{t.cell}
 Let $T$ be a standard $D$-Young tableau and let $c=cont(T)$. Then for any $\pi \in D_n$,
$$\pi \in K^c \Longleftrightarrow \textrm{ The tableau } T^{\pi^{-1}} \text{
is a standard $D$-Young tableau}.$$
\end{thm}

\section{Standard D-Young tableaux: construction and
enumeration.}

 We deal here only with tableaux which have at most two boxes on
the zero content diagonal. The reason for it is that tableaux
without boxes with zero content lead to representations induced to
$D_n$ from $S_n$, while tableaux with one or two boxes of zero
content give rise to minimal AY representations which are
subrepresentations of these induced representations.

\subsection{Tableaux without zero-content boxes }

\bde Let $m\in \mathbb N$. Denote by $T[\lambda,m,+]$ the set
which consists of all standard D-Young tableaux whose smallest
nonnegative content is equal to $m$, whose boxes with positive
contents form the shape $\lambda$ and the number of negative
entries in the boxes with positive contents is even. The set
$T[\lambda,m,-]$ is defined similarly to $T[\lambda,m,+]$, just
the number of negative entries in the boxes with positive contents
is odd. \ede

Further, when we consider both $T[\lambda,m,+]$ and $T[\lambda,m,-]$
we sometimes write $T[\lambda,m,\pm]$ for brevity.

Obviously from the definition, tableaux from sets
$T[\lambda,m,\pm]$ have no zero-content boxes.

Recall that if the number of negative entries in the boxes of $T$
with positive contents is even, then it is also even for $T^{\pi}$
with any $\pi\in D_n$. For example, if we start with
$$
T=\,\,\young(:::::123,\mthree\mtwo\mone)\qquad cont(T)=(2,3,4)
$$
we obtain the 4-element set $T[(3),2,+]$.

\bpr\label{inddim} Let $\lambda$ be a shape (straight or skew)
with $n$ boxes and $m\in\mathbb N$. Then
$$\# T[\lambda,m,+]=\# T[\lambda,m,-]=2^{n-1}f^{\lambda}$$
\epr
\begin{proof}  We are dealing now with
tableaux that have no box on zero-content diagonal. This means
that the standardness (i.e. increasing of the entries in rows and
columns) of the sub tableau with positive contents is not affected
at all by the part with negative contents. In other words, when
there is no box with content zero, then one can build a standard
D-Young tableau just filling the shape $\lambda$ (which is the
shape of $n$ boxes with positive contents) in the standard way
(i.e. increasing in rows and columns) by the numbers
$\{x_1,...,x_n\}$ such that $\{|x_1|,\cdots ,|x_n|\}=\{1,2,\cdots,
n\}$. For each choice of such numbers $\{x_1,...,x_n\}$ there are
$f^{\lambda}$ ways to arrange them in the shape $\lambda$
increasing in rows and columns. For the tableaux of the set
$T[\lambda,m,+]$ the number of negative entries in the shape
$\lambda$ has to be even i.e. it must be 0,2,4,... and for
$T[\lambda,m,-]$ this number must be 1,3,5,... Therefore
$$
\# T[\lambda,m,+]=\left({n\choose 0}+{n\choose 2}+{n\choose
4}+\cdots\right)f^{\lambda}=2^{n-1}f^{\lambda}
$$
and
$$
\# T[\lambda,m,-]=\left({n\choose 1}+{n\choose 3}+{n\choose
5}+\cdots\right)f^{\lambda}=2^{n-1}f^{\lambda}\,\,.
$$

\end{proof}

\subsection{Tableaux having one or two boxes with zero content}
Now consider tableaux with at most two boxes on the
diagonal of zero content. Such tableau can be divided (sometimes
not uniquely) into two sub tableaux: one contains $n$ boxes with
nonnegative contents and another (which is a reflection of the
first one because $c_{-i}=-c_i$) contains the boxes with non
positive contents.

\subsubsection{Tableaux having one box with zero content}
First consider the case in which the tableaux have one box with
zero content. The condition $c_{-i}=-c_i$ implies that this unique
zero-content box must be occupied by two numbers $\pm i$. For
example,
$$
\young(:\mthree1,:\pmtwo,\mone3)
$$
Here the shape (but not a sub tableau!) of nonnegative contents
can be defined uniquely and is
$$
\lambda=\young(\,\,,\,)=(2,1)
$$
\bde We define the set $T[\lambda,\boxtimes]$ to be the set of all
standard D-Young tableaux which have one box with zero content
(always occupied by $\pm i$ for some $i$) and shape of the boxes
with nonnegative contents is $\lambda$. \ede

For example, the set $T[(2,1),\boxtimes]$ is
$$
\young(:\mthree\mtwo,:\pmone,23)\,\,\,;
\,\,\young(:\mthree1,:\pmtwo,\mone3)\,\,\,;
\,\,\young(:\mthree\mone,:\pmtwo,13)\,\,\,;
\,\,\young(:\mthree2,:\pmone,\mtwo3)\,\,\,;
\,\,\young(:\mtwo3,:\pmone,\mthree2)
$$
while the set $T[(4),\boxtimes]$ consists of only one tableau $
\young(\mfour\mthree\mtwo\pmone234)$.

\subsubsection{Tableaux having two boxes with zero content}

Now suppose we have a tableau with two different boxes on the zero
content diagonal.

\bde\label{verhor} Denote by $T[\lambda,\cdot | \cdot,+]$
($T[\lambda,\div,+]$) the set of standard D-Young tableaux which
have exactly two boxes on the zero content diagonal  and which can
be divided by the \textbf{vertical} (\textbf{horizontal}) straight
line into two parts (with $n$ boxes in each part)-- one consists
of the boxes with nonnegative contents and another of non positive
contents -- and $\lambda$ is the shape (straight or skew) of boxes
with nonnegative contents after this separation and the number of
negative entries in the boxes of $\lambda$ is even. When the
number of negative entries in the boxes of $\lambda$ is odd we
have $T[\lambda,\cdot |\cdot,-]$ ($T[\lambda,\div,-]$).
 \ede
Consider, for example, the following tableau:
$$
\young(:\mthree12,\mtwo\mone3)
$$
The diagonal of zero contents passes through the number 3 and $-3$
and this tableau can be divided as required above by the vertical
line:
$$
\young(:\mthree12,\mtwo\mone3)\,\,\,\,\longrightarrow\,\,\,\,\,\young(:\mthree,\mtwo\mone)\qquad\young(12,3)
$$
Here $\lambda=\young(\,\,,\,)=(2,1)$ and the number of negative
entries in $\lambda$ is equal to zero, thus even. The complete set
$T[(2,1),\cdot | \cdot,+]$ is
$$
\young(:\mthree12,\mtwo\mone3)\,\,\,;
\,\,\,\young(:\mtwo13,\mthree\mone2)\,\,\,; \,\,\,
\young(:\mthree\mtwo\mone,123)
$$
while the set $T[(2,1),\cdot | \cdot,-]$ is
$$
\young(:\mthree\mone2,\mtwo13)\,\,\,;
\,\,\,\young(:\mtwo\mone3,\mthree12)\,\,\,; \,\,\,
\young(:\mthree\mtwo1,\mone23)
$$

Not every standard D-young tableau can be divided by a vertical
line into two halves of nonnegative and non positive contents. For
example, if we divide by the vertical line the following tableau
$$
\young(\mfour\mthree,\mtwo1,\mone2,34)\,\,\,\,\longrightarrow\,\,\,\,\,\young(\mfour,\mtwo,\mone,3)
\qquad\young(\mthree,1,2,4)
$$
each half has positive and negative contents together.

The following tableau can be divided as required above by the
horizontal line:
$$
\young(:\mthree,\mtwo\mone,12,3)\,\,\,\,\longrightarrow
\,\,\,\,\,\begin{matrix}\young(:\mthree,\mtwo\mone)\\
\\ \young(12,3)\end{matrix}
$$
Here the skew shape $\lambda=\young(:\,,\,\,)=(2,2)/(1)$ and the
complete set is
$$
T[(2,2)/(1),\div,-]=\left\{\,\,\,\young(:\mthree,\mtwo\mone,12,3)\,\,\,\,;
\,\,\,\young(:\mtwo,\mthree\mone,13,2)\,\,\,\,;
 \,\,\,\young(:1,\mthree2,\mtwo3,\mone)\,\,\,\,\right\}
 $$

The structure of the D-tableaux with two boxes on the zero
diagonal which belong to both $T[\lambda,\div,+]$ and
$T[\tilde{\lambda},\cdot |\cdot,-]$ is described in the following
observation (this observation will not be used later):
 \bobs Let $\lambda=\mu/\nu$ for some partitions
$\mu=(\mu_1,\ldots,\mu_h)\vdash n+k$ where $h=\mu'_1$ and
$\nu=(\nu_1,\ldots,\nu_{\ell})\vdash k$ where $\ell=\nu'_1$, and
$h>\ell$ and suppose that $h-\ell=1$ and $\mu_h\geqslant 2$. Then
$$
T[\lambda,\div,+]=T[\tilde{\lambda},\cdot
|\cdot,-]\quad\textrm{and}\quad
T[\lambda,\div,-]=T[\tilde{\lambda},\cdot |\cdot,+]
$$
where $\tilde{\lambda}=\tilde{\mu}/\tilde{\nu}$ for
$\tilde{\mu}=(\tilde{\mu}_1,\ldots,\tilde{\mu}_h,\tilde{\mu}_{h+1})$
with $\tilde{\mu}_1=\mu_1-1$ , $\tilde{\mu}_2=\mu_2-1$ , ...,
$\tilde{\mu}_h=\mu_h-1$ , $\tilde{\mu}_{h+1}=1$ and
$\tilde{\nu}=(\tilde{\nu}_1,\ldots,\tilde{\nu}_{\ell})$ with
$\tilde{\nu}_1=\nu_1-1$ , $\tilde{\nu}_2=\nu_2-1$ ,...,
$\tilde{\nu}_{\ell}=\nu_{\ell}-1$. Notice that the number of parts
in the partition $\tilde{\nu}$ can be less than $\ell$. Moreover,
if $\nu=1^k$, then $\tilde{\nu}$ is empty. \eobs

We finish this subsection with a simple remark which easily
follows from the definition of $T[\lambda,\cdot |\cdot,\pm]$ and
$T[\lambda,\div,\pm]$:
\begin{rem}
For each standard D-Young with exactly two boxes on the zero
diagonal, there exists a shape $\lambda$ such that this tableau
belongs to $T[\lambda,\cdot |\cdot,\pm]$ or $T[\lambda,\div,\pm]$.
\end{rem}
\subsubsection{Constructing D-Young tableaux from a shape with $n$
boxes.}

We continue to study sets $T[\lambda,\cdot |\cdot,\pm]$,
$T[\lambda,\div,\pm]$, $T[\lambda,\boxtimes]$ and
$T[\lambda,m,\pm]$ in more detail.
 Lemmas~\ref{T0} and~\ref{Tm} which are proved
in this subsection are needed for the proof of the main result --
Theorem~\ref{decomp}.

 Suppose we are given a shape (straight or skew) $\lambda$
with $n$ boxes. We shall address the following question: what sets
of standard D-Young tableaux may be constructed by putting the
lower left box of $\lambda$ on the diagonal of zero content?

Obviously it is always (i.e. for any $\lambda$) possible to
construct the set $T[\lambda,\boxtimes]$. The conditions which
provide the existence of $T[\lambda,\cdot |\cdot,\pm]$ and
$T[\lambda,\div,\pm]$ are in given in the following simple remark.
\begin{rem}\label{T0sets} If $h-\ell\geqslant 2$ and $\mu_h=1$, then the sets
$T[\lambda,\cdot |\cdot,\pm]$ are defined but the sets
$T[\lambda,\div,\pm]$ are not defined and we assume
$T[\lambda,\div,\pm]=\varnothing$. If $h-\ell=1$ and
$\mu_h\geqslant 2$, then the sets $T[\lambda,\div,\pm]$ are
defined but the sets $T[\lambda,\cdot |\cdot,\pm]$ don't exist and
we assume $T[\lambda,\cdot |\cdot,\pm]=\varnothing$. If
$h-\ell\geqslant 2$ and $\mu_h\geq 2$, then both
$T[\lambda,\div,\pm]$ and $T[\lambda,\cdot |\cdot,\pm]$ are
available. If $\mu_h=1$ and $h-\ell=1$, then the only thing we can
construct putting the lower left box of $\lambda$ on the diagonal
of zero content is the set $T[\lambda,\boxtimes]$ while the sets
$T[\lambda,\cdot |\cdot,\pm]$ and $T[\lambda,\div,\pm]$ are empty.
\end{rem}
\bde\label{zerotab} Denote
$$
T[\lambda,0,+]:=T[\lambda,\cdot |\cdot,+]\bigcup
T[\lambda,\div,+]\bigcup T[\lambda,\boxtimes]
$$
and
$$
T[\lambda,0,-]:=T[\lambda,\cdot |\cdot,-]\bigcup
T[\lambda,\div,-]\bigcup T[\lambda,\boxtimes]\,\,\,.
$$
\ede

The following theorem describes a bijection between
$T[\lambda,0,\pm]$ and $T[\lambda,m,\pm]$ with $m=1,2,3,...$
\begin{thm}\label{biject}
 Let $\lambda$ be a skew or straight shape with $n$ boxes. Then for any $m\in\mathbb N$ there
 exists a natural bijection:
$$
T[\lambda,0,+]=T[\lambda,\boxtimes]\bigcup
T[\lambda,\div,+]\bigcup T[\lambda,\cdot
|\cdot,+]\,\,\,\longleftrightarrow\,\,\, T[\lambda,m,+]
$$
and
$$
T[\lambda,0,-]=T[\lambda,\boxtimes]\bigcup
T[\lambda,\div,-]\bigcup T[\lambda,\cdot
|\cdot,-]\,\,\,\longleftrightarrow \,\,\,T[\lambda,m,-]\,\,.
$$
\end{thm}
Theorem~\ref{biject} will be proved using Lemmas~\ref{T0}
and~\ref{Tm}. It should be noted that Lemma \ref{T0} will be used
later, when analyzing the associated representations.
\blem\label{T0} Let $\lambda$ be a straight or a skew shape with
$n$ boxes, i.e. $\lambda =\mu/\nu$ for some partitions
$\mu=(\mu_1,\ldots,\mu_h)\vdash n+k$ where $h=\mu'_1$ and
$\nu=(\nu_1,\ldots,\nu_{\ell})\vdash k$ where $\ell=\nu'_1$, and
$h>\ell$, and let $T$ be a standard Young tableau of shape
$\lambda$ filled by entries $x_1$, $x_2$,..., $x_n$ such that
$\{|x_1|,|x_2|,...,|x_n|\}=\{1,2,...,n\}$. If the number of
negative entries in $T$ is even (odd), then there exists a unique
standard D-Young tableau $\lceil T\rfloor_0\in T[\lambda,0,+]$
$(\lceil T\rfloor_0\in T[\lambda,0,-]$, respectively$)$ such that
$T$ is its half of shape $\lambda$ with nonnegative contents.
\elem
\begin{proof} (1) Let $h-\ell\geqslant 2$ and $\mu_h\geqslant 2$. Then
the lower left corner of $T$ has at least two boxes in both
directions:
$$
T=\young(:**,a**,bc)
$$
where $a<b<c$.

Denote
$$
X=\young(::::**,:\mc\mb a**,**\ma
bc,**)\,,\,Y=\young(::**,:a**,:bc,:\mc\mb,**\ma,**)\,,\,Z=\young(:::**,::a**,:\mc\plmb
c,**\ma,**)\,\,.
$$

 If $|a|<|b|<|c|$, then $X$ is a standard D-Young tableau while $Y$ and $Z$ are obviously non standard and we put
$\lceil T\rfloor_0=X\in\,\, T[\lambda,\cdot |\cdot,+]$.

If $|a|>|b|>|c|$, then $Y$ is standard while $X$ and $Z$ aren't
and so we put $\lceil T\rfloor_0=Y\in\,\, T[\lambda,\div,+]$.

 Finally, if $|a|>|b|<|c|$, then $Z$ is standard while $X$ and $Y$
 aren't and
so we put $\lceil T\rfloor_0=Z\in\,\, T[\lambda,\boxtimes]$.

 The case $|a|<|b|>|c|$ is impossible because $a<b<c$.

Consider the examples. Take
$$
T=\,\,\young(\mtwo\mone,34)\,\,\,\textrm{with its
reflection}\,\,\, \young(\mfour\mthree,12)\,.
$$
Then $\lceil
T\rfloor_0=\young(\mfour\mthree\mtwo\mone,1234)\,\,\in
T[(2,2),\cdot |\cdot,+]$ while the tableaux
$$
\young(\mtwo\mone,34,\mfour\mthree,12)
\,\,\,\,\,\,\,\,\,\,\,\textrm{and}\,\,\,\,\,\,\,\,\,\young(:\mtwo\mone,\mfour\pmthree4,12)
\,\,\,\,\,\,\,\textrm{are obviously not standard.}
$$
 Two more examples:
$$
T=\young(\mfour\mthree,\mtwo\mone)
\,\,\,\,\Longrightarrow\,\,\,\lceil
T\rfloor_0=\young(\mfour\mthree,\mtwo\mone,12,34)\,\,\,\in
\,\,T[(2^2),\div,+]
$$
and
$$
T=\young(\mfour2,\mone3)\,\,\,\,\Longrightarrow \,\,\,\lceil
T\rfloor_0=\young(:\mfour2,\mthree\pmone3,\mtwo4)\,\,\,\in
\,\,T[(2^2),\boxtimes]\,\,.
$$
\\

\noindent(2) If $h-\ell=1$ and $\mu_h\geqslant 2$, then  the box
$\young(a)$ is absent in $T$ and the tableau $X$ doesn't exist,
$T[\lambda,\cdot|\cdot,+]=\varnothing$. If $h-\ell\geqslant 2$ and
$\mu_h=1$, then the box $\young(c)$ is absent in $T$ and the
tableau $Y$ doesn't exist, $T[\lambda,\div,+]=\varnothing$. The
above item (1) of our proof can be repeated with slight changes:
we have to erase $a$ or $c$, respectively, from all the places
where they appear. If $h-\ell=1$ and $\mu_h=1$, then both
$\young(a)$ and $\young(c)$ are absent, both of the sets
$T[\lambda,\cdot |\cdot,+]$ and $T[\lambda,\div,+]$ are empty, the
lower left corner of $T$ consists of a single box and $ \lceil
T\rfloor_0\in T[\lambda,\boxtimes]$.

Thus, for any standard Young tableau $T$ of an arbitrary shape
$\lambda$ with $n$ boxes filled by entries $x_1$, $x_2$,..., $x_n$
  such that
$\{|x_1|,|x_2|,...,|x_n|\}=\{1,2,...,n\}$ and the number of
minuses in $T$ is even we have uniquely defined a standard D-Young
tableau $\lceil T\rfloor_0\in T[\lambda,0,+]$ such that $T$ is its
half of shape $\lambda$ with nonnegative contents, as claimed.
When the number of negative entries in $T$ is odd, we can repeat
with appropriate slight changes all the arguments of the above
proof to define $\lceil T\rfloor_0\in T[\lambda,0,-]$.
\end{proof}

\blem\label{Tm} Let $\lambda$ be a straight or a skew shape with
$n$ boxes and let $T$ be a standard Young tableau of shape
$\lambda$ filled by entries $x_1$, $x_2$,..., $x_n$
 such that
$\{|x_1|,|x_2|,...,|x_n|\}=\{1,2,...,n\}$. If the number of
negative entries in $T$ is even (odd), then for any $m\in\mathbb
N$ there exists a unique standard D-Young tableau $\lceil
T\rfloor_m\in T[\lambda,m,+]$ $(\lceil T\rfloor_m\in
T[\lambda,m,-]$, respectively$)$ such that $T$ is its half of
shape $\lambda$ with positive contents. \elem
\begin{proof}
Take $T$ and put it in such a way that the box at its lower left
corner has the content $m$. Then put $T$'s reflection in such a
way that its upper right corner box has the content $-m$ and the
whole picture is a skew shape, i.e. $T$ is strictly to the
northeast of its reflection. Clearly we obtain a D-Young tableau.
Since there is no box on the zero content diagonal, the halves of
positive and negative contents don't touch each other and
therefore the obtained D-Young tableau is standard and we denote
it as $\lceil T\rfloor_m$. Example:
$$\lambda=(2)\vdash 2
\,\,,\,\,T=\young(\mone2)\,\,,m=1\,\,\Longrightarrow\lceil
T\rfloor_1=\young(::\mone2,\mtwo1)\,\,\in T[(2),1,-]\nonumber
$$
\end{proof}
From Lemmas~\ref{T0} and~\ref{Tm}, we can easily conclude the
following

\noindent\textbf{Proof of Theorem \ref{biject}.} For any standard
D-Young tableau $X\in T[\lambda,0,\pm]$ or $T[\lambda,m,\pm]$, we
have $X=\lceil T\rfloor_0$ or $\lceil T\rfloor_m$ where $T$ is the
half of $n$ boxes with nonnegative (positive) contents of shape
$\lambda$ in $X$, and our claim now follows.\qed

\section{Several types of DAY representations.}\label{dnconst}

Given a standard D-Young tableau $T$ with its content vector
$$c(T)=(c_1(T),c_2(T)...,c_n(T))\,.$$
As discussed in previous sections, $c(T)$ is a generic vector (see
Theorem~\ref{tablestoreps}) and thus it gives rise to a minimal AY
cell in $D_n$ (see Theorem~\ref{t.main11}). Denote
\begin{align*}
f(T)&=(f_0(T),f_1(T),...,f_{n-1}(T))\\
&=\Delta c(T)=(c_1+c_2,c_2-c_1,c_3-c_2,...,c_n-c_{n-1})\,.
\end{align*}

 The corresponding
minimal AY representation $\rho$, which we called a $DAY$
representation (see Definition~\ref{DAY}) acts on the space
spanned by all the standard D-Young tableaux obtained from $T$ by
the natural action of $D_n$. The action of representation matrices
of generators is defined by
$$\rho_{s_i}(T)=\frac{1}{f_i(T)}\,\,T+\sqrt{1-\frac{1}{f_i(T)^2}}\,\,T^{s_i}
\,\,\,\,\textrm{for}\,\,i=0,1,2,...,n-1\,.\eqno(*)
$$
Notice (see Theorem~\ref{t.cell}) that $T^{s_i}$ is not standard
if and only if
$$
\sqrt{1-\frac{1}{f_i(T)^2}}\left(=\frac{\sqrt{f_i(T)^2-1}}{|f_i(T)|}\right)=0
\Longleftrightarrow f_i=\pm1\,\,.
$$
\bde For $m\in \{0\}\bigcup\mathbb N$ we denote by
$\rho^{\lambda,m,+}$, $\rho^{\lambda,m,-}$, $\rho^{\lambda,\cdot
|\cdot,+}$, $\rho^{\lambda,\cdot |\cdot,-}$,
$\rho^{\lambda,\div,+}$, $\rho^{\lambda,\div,-}$,
$\rho^{\lambda,\boxtimes}$ DAY representations defined above by
$(*)$ on spaces spanned by sets $T[\lambda,m,+]$,
$T[\lambda,m,-]$, $T[\lambda,\cdot |\cdot,+]$, $T[\lambda,\cdot
|\cdot,-]$, $T[\lambda,\div,+]$, $T[\lambda,\div,-]$,
$T[\lambda,\boxtimes]$, respectively.

\ede

\begin{rem}\label{zerorep}
The sets $T[\lambda,\cdot |\cdot,+]$ , $T[\lambda,\div,+]$ ,
$T[\lambda,\boxtimes]$ are pairwise disjoint because they consist
of tableaux with different shapes and therefore
$$
\rho^{\lambda,0,+}=\rho^{\lambda,\cdot
|\cdot,+}\oplus\rho^{\lambda,\div,+}\oplus\rho^{\lambda,\boxtimes}
$$
and similarly
$$
\rho^{\lambda,0,-}=\rho^{\lambda,\cdot
|\cdot,-}\oplus\rho^{\lambda,\div,-}\oplus\rho^{\lambda,\boxtimes}\,\,\,.
$$
\end{rem}
When we write $\rho^{\lambda,m,\pm}$, $\rho^{\lambda,\cdot
|\cdot,\pm}$ or $\rho^{\lambda,\div,\pm}$ we mean
``representations $\rho^{\lambda,m,+}$ and $\rho^{\lambda,m,-}$"
or ``representations $\rho^{\lambda,\cdot |\cdot,+}$ and
$\rho^{\lambda,\cdot |\cdot,-}$" or ``representations
$\rho^{\lambda,\div,+}$ and $\rho^{\lambda,\div,-}$",
respectively.

\subsection{Representations $\rho^{\lambda,m,\pm}$ are induced
from $S_n$ to $D_n$.}\label{indSntoDn}

There are two relevant embeddings of $S_n$ in $D_n$. Denote them
$$
S_n^1=\langle s_1,s_2,\ldots,s_{n-1}\rangle
$$
and
$$
S_n^0=\langle s_0,s_2,\ldots,s_{n-1}\rangle\quad.
$$
When $n$ is odd $S_n^1$ and $S_n^0$ are conjugate subgroups in
$D_n$. When $n$ is even they are not conjugate.

\bigskip
There is the following known fact:
\begin{fac}\label{cosets1}
The cosets are $D_n/S_n^1=\{\sigma_1 S_n^1,\sigma_2
S_n^1,\ldots,\sigma_{2^{n-1}} S_n^1\}$ where\\
$\sigma_1=e=(1\,\,2\,\,3\,\,4\,\ldots\, n)\,\,, $\\
$\sigma_2=s_0=(-2\,\,-1\,\,3\,\,4\,\ldots\, n)\,\,,$\\
$\sigma_3=s_2s_0=(-3\,\,-1\,\,2\,\,4\,\ldots\,n)\,\,,$\\
$\sigma_4=s_1s_2s_0=(-3\,\,-2\,\,1\,\,4\,\ldots\,n)\,\,,$\\
$\sigma_5=s_3s_2s_0=(-4\,\,-1\,\,2\,\,3\,\ldots\,n)\,\,,$\\
$\sigma_6=s_1s_3s_2s_0=s_3s_1s_2s_0=(-4\,\,-2\,\,1\,\,3\,\ldots\,n)\,\,,$\\
$\sigma_7=s_2s_1s_3s_2s_0=(-4\,\,-3\,\,1\,\,2\,\ldots\,n)\,\,,$\\
$\sigma_8=s_0s_2s_1s_3s_2s_0=(-4\,\,-3\,\,-2\,\,-1\,\ldots\,n)\,\,,$\\
$\cdots\cdots\cdots\cdots\cdots\cdots\cdots\cdots\cdots\cdots\cdots$\\
$\sigma_{2^{n-1}}=(-n\,\,-(n-1)\,\,-(n-2)\,\ldots\,-2\,\,-1)$ if
$n$ is even and\\
$\sigma_{2^{n-1}}=(-n\,\,-(n-1)\,\,-(n-2)\,\ldots\,-2\,\,1)$ if
$n$ is odd.\\
 These $\sigma_i$ are coset representatives of minimal length.
\end{fac}

\begin{fac}\label{cosetminus}
To obtain the minimal length representatives of cosets of
$D_n/S_n^0$ which we shall denote as $\tilde{\sigma}_i$ we just
have to change the roles of $s_1$ and $s_0$ in each $\sigma_i$. It
is convenient to index them by the increasing sequences  with odd
number of minuses which are obtained by the action of each $\tilde{\sigma}_i$
on the sequence $(-1\,\,2\,\,3\,\ldots\, n)$:\\
$\tilde{\sigma}_1=e \longleftrightarrow(-1\,\,2\,\,3\,\,4\,\ldots\, n)\,\,, $\\
$\tilde{\sigma}_2=s_1\longleftrightarrow(-2\,\,1\,\,3\,\,4\,\ldots\, n)\,\,,$\\
$\tilde{\sigma}_3=s_2s_1\longleftrightarrow(-3\,\,1\,\,2\,\,4\,\ldots\,n)\,\,,$\\
$ \tilde{\sigma}_4=s_0s_2s_1\longleftrightarrow(-3\,\,-2\,\,-1\,\,4\,\ldots\,n)\,\,,$\\
$ \tilde{\sigma}_5=s_3s_2s_1\longleftrightarrow(-4\,\,1\,\,2\,\,3\,\ldots\,n)\,\,,$\\
$\cdots\cdots\cdots\cdots\cdots\cdots\cdots\cdots\cdots\cdots\cdots$
\end{fac}
Now we prove that for any $m\in \mathbb N$ our minimal AY
representations $\rho^{\lambda,m,+}$ are induced into $D_n$ from
$S_n^1$ while $\rho^{\lambda,m,-}$ are induced into $D_n$ from
$S_n^0$.

 Let $\lambda$ be a straight or skew shape with $n$ boxes and let
$S^{\lambda}$ denote the representation of $S_n$ associated with
$\lambda$ via classical Young orthogonal form (see
Theorems~\ref{t.yof} and~\ref{t.yof-skew}) . (If $\lambda$ is a
straight shape, then $S^{\lambda}$ is an irreducible Specht module
and if $\lambda$ is a skew shape, then $S^{\lambda}$ is a skew
Specht module.)

\begin{thm}\label{ind}
 For any natural $m=1,2,3,\ldots$
 \begin{itemize}
 \item[(1)]   $\quad\quad\rho^{\lambda,m,+}\cong
 S^{\lambda}\uparrow_{S_n^1}^{D_n}$

\item[(2)]   $\quad\quad\rho^{\lambda,m,-}\cong
 S^{\lambda}\uparrow_{S_n^0}^{D_n}$
\end{itemize}
\end{thm}

\begin{proof}
We give the proof of (1) in detail and, then briefly explain how
the proof of (2) follows.

Recall that the representation space of $\rho^{\lambda,m,+}$ is
spanned by the set of standard D-Young tableaux $T[\lambda,m,+]$.
As we already know (see Lemma~\ref{Tm}) each  standard tableau $T$
which is a filling of $\lambda$ (which is the shape of boxes with
positive contents) by $n$ numbers from the set
$\{\pm1,\pm2,\ldots,\pm n\}$ in such a way that from each pair
$\pm i$ only one number (i.e. or $i$ or $-i$) appears in the
tableau and the number of negatives in the tableau is even leads
to the unique standard D-Young tableau $\lceil T\rfloor_m\in
T[\lambda,m,+]$ for which $T$ is its sub tableau of shape
$\lambda$ of boxes with positive contents. So we can work with
standard tableaux $T$ of shape $\lambda$ disregarding the parts
with negative contents and from now on in this subsection we
identify a standard D-Young tableau $\lceil T\rfloor_m\in
T[\lambda,m,+]$ with $T$.

From now on in this subsection we denote for brevity
$\rho=\rho^{\lambda,m,+}$.

Fix some standard Young tableau $T$ of shape $\lambda$ filled with
numbers 1,2,3,...,$n$. Notice that for any other standard Young
tableau $\widetilde{T}$ of shape $\lambda$ filled with numbers
1,2,3,...,$n$ there exists $\pi\in S_n^1$ such that
$\widetilde{T}=T^{\pi}$ (see proof of Theorem \ref{t.cell} for more
details). Obviously the subspace of $Span\{T[\lambda,m,+]\}$ spanned
by all (standard D-Young tableaux $\lceil T\rfloor_m$ whose halves
of positive contents are) standard Young tableaux of shape $\lambda$
filled with numbers 1,2,3,...,$n$ is invariant under the action of
$\rho_{s_i}$ for $0<i=1,2,\ldots,n-1$ (i.e. is invariant under the
action of $S_n^1$) and the corresponding representation is exactly
the representation $S^{\lambda}$ of $S_n^1$.

The following obvious observation is important:

\bobs\label{ind1} Let
$SYT(\lambda)=\{T_1,T_2,\ldots,T_{f^{\lambda}}\}$ be the set of all
standard Young tableaux of shape $\lambda$ filled with numbers
1,2,3,...,$n$. Then all the tableaux $T_j^{\sigma_i}$ for
$j=1,2,\ldots,f^{\lambda}$ and $i=1,2,\ldots,2^{n-1}$ with
$\sigma_i$ defined in Fact \ref{cosets1} are also standard, i.e.
their entries are increasing in rows and columns, and
$T_j^{\sigma_i}\neq T_{j'}^{\sigma_i}$ for $j\neq j'$. Moreover, it
is easy to see by induction on $\ell(\sigma_{i})$ (the length of
$\sigma_{i}$) that there exist
 signed permutations
$\omega\in D_n$ and coefficients $\alpha^m_{j,i,\omega}\neq 0$,
$\alpha^m_{j,i,\sigma_{i} }\neq 0$ such that tableaux
$T_{j}^{\omega}$ are standard and
$$
\rho_{\sigma_{i}}(T_{j})=\sum\alpha^m_{j,i,\omega}
T_{j}^{\omega}+\alpha^m_{j,i,\sigma_{i}}T_{j}^{\sigma_{i}}\quad\textrm{with}\quad
\ell(\sigma_{i})>\ell(\omega)\quad\textrm{for any}\quad\omega\,.
$$
\eobs

 In order to prove the first statement of Theorem \ref{ind} we must show that
$$
Span\{T[\lambda,m,+]\}=\bigoplus_{i=1}^{2^{n-1}}\rho_{\sigma_i}(Span\{SYT(\lambda)\})\,\,
$$
where
$$\rho_{\sigma_i}(Span\{SYT(\lambda)\})
=Span\{\rho_{\sigma_i}(T_1),\rho_{\sigma_i}(T_2),
\rho_{\sigma_i}(T_3),\ldots,\rho_{\sigma_i}(T_{f^{\lambda}})\}\,\,.
$$
The number $m$ affects only on the certain values of coefficients
$\alpha^m_{j,i,\omega}$ but all the arguments of the proof remain
the
same for any $m\in \mathbb N$.\\

Now we explain why the subspaces
$\rho_{\sigma_i}(Span\{SYT(\lambda)\})$ and
$$\rho_{\sigma_{i'}}(Span\{SYT(\lambda)\})
=Span\{\rho_{\sigma_{i'}}(T_1),\rho_{\sigma_{i'}}(T_2),
\rho_{\sigma_{i'}}(T_3),\ldots,\rho_{\sigma_{i'}}(T_{f^{\lambda}})\}$$
have trivial intersection:

\blem\label{trivintersect} In current notations
$$
\rho_{\sigma_i}(Span\{SYT(\lambda)\})\bigcap\rho_{\sigma_{i'}}(Span\{SYT(\lambda)\})=\{0\}\,\,\,\,\textrm{for
$i\neq i'$.}
$$
\elem

 The proof of this lemma will
follow from two claims. \bcl\label{cl1} $T_j^{\sigma_i}$ cannot
appear in the decomposition of $\rho_{\sigma_{i'}}(T_{j'})$ for
$i'<i$ and\linebreak any $j'$.\ecl

\noindent \textbf{Proof of Claim \ref{cl1}}
 Let the coset
representatives $\sigma_i$ be arranged by\linebreak increasing of
length, i.e. $\ell(\sigma_i)\geqslant \ell(\sigma_{i'})$ when
$i>i'$.

  Recall that $T_j,T_{j'}\in SYT(\lambda)$ and hence
$T_{j'}=\left(T_j\right)^{\pi}$ for some $\pi\in S_n^1$ (if
$j'=j$, then $\pi=e$) and assume to the contrary that
$T_j^{\sigma_i}$ appears in the decomposition of
$\rho_{\sigma_{i'}}(T_{j'})$ for some $i'<i$. Then by
Observation~\ref{ind1}, we must have
$T_j^{\sigma_i}=T_{j'}^{\sigma_{i'}}$ or
$T_j^{\sigma_i}=T_{j'}^\omega$ for some certain $\omega\in D_n$
with $\ell(\sigma_{i})\geqslant\ell(\sigma_{i'})>\ell(\omega)$.
The equality $T_j^{\sigma_i}=T_{j'}^{\sigma_{i'}}$ is impossible
because
$T_j^{\sigma_i}=T_{j'}^{\sigma_{i'}}=\left(\left(T_j\right)^{\pi}\right)^{\sigma_{i'}}=
\left(T_j\right)^{\sigma_{i'}\pi}$ for some $\pi\in S_n^1$ means
that ${\sigma_i}$ and $\sigma_{i'}$ belong to the same coset which
is a contradiction. The equality
$T_j^{\sigma_i}=T_{j'}^\omega=T_j^{\omega\pi}$ with $\pi\in
S_n^1$, $\omega\in D_n$ and $\ell({\sigma_i})>\ell(\omega)$ is
also impossible because it means that $\sigma_i=\omega\pi$ and
hence $\sigma_i$ is not a coset representative of minimal length:
$\omega=\sigma_i\pi^{-1}\in\sigma_i S_n^1$ and has shorter length
than $\sigma_i$. This contradiction completes the proofs that
$T_j^{\sigma_i}$ can not appear in the decomposition of
$\rho_{\sigma_{i'}}(T_{j'})$ for $i'<i$ and any $j'$ and thus
Claim~\ref{cl1} is proved.\qed

\bcl\label{cl2} Let $\{u_1,u_2,...,u_k\}$ and
$\{v_1,v_2,...,v_k\}$ be two linearly independent subsets of a
vector space $U$ such that $Span\{u_i\}_{i=1}^{k}\bigcap
Span\{v_i\}_{i=1}^{k}=\{0\}$. Then
$$
Span\{u_1,u_2,...,u_k\}\bigcap
Span\{u_1+v_1,u_2+v_2,...,u_k+v_k\}=\{0\}\,\,.
$$
\ecl \noindent The proof of Claim~\ref{cl2} is a simple exercise
in linear algebra.\qed

 \noindent \textbf{Proof of Lemma \ref{trivintersect}} Combine
Observation~\ref{ind1} with Claims~\ref{cl1} and~\ref{cl2}.\qed

 The subspace
$\rho_{\sigma_1}(Span\{SYT(\lambda)\})=Span\{SYT(\lambda)\}$ is
invariant under the action of $\rho_{s_1}$ , $\rho_{s_2}$ , ... ,
$\rho_{s_{n-1}}$ (it is an $S_n^1$-representation $S^{\lambda}$).
By the construction, the $\rho$-action of $D_n$ on
$Span\{T[\lambda,m,+]\}$ permutes the subspaces
$\rho_{\sigma_i}(Span\{SYT(\lambda)\})$ and because these subspace
have trivial intersection their sum is direct. Also clear that for
each $i$ we have $\dim
\rho_{\sigma_i}(Span\{SYT(\lambda)\})=f^{\lambda}$ and so
$$
\bigoplus_{i=1}^{2^{n-1}}\rho_{\sigma_i}(Span\{SYT(\lambda)\})=Span\{T[\lambda,m,+]\}
$$
All what we said above means that
$\quad\quad\rho^{\lambda,m,+}\cong
 S^{\lambda}\uparrow_{S_n^1}^{D_n}$ by the definition of induced
 representation and thus the proof of the statement (1) of Theorem \ref{ind} is completed.

 The proof of the statement (2) of Theorem \ref{ind} is very similar to above.
 Instead of $SYT(\lambda)$
 filled by numbers 1,2,...,$n$ take the set
 $\overline{SYT(\lambda)}$ of standard tableaux of shape $\lambda$
 filled by numbers $-1$,2,...,$n$. Easy to see that
 the subspace $Span\{\overline{SYT(\lambda)}\}$ is invariant under the
 action of $\rho_{s_0}$ , $\rho_{s_2}$ , $\rho_{s_3}$
 ,..., $\rho_{s_{n-1}}$ and gives an $S_n^0$-representation
 $S^{\lambda}$. We list the coset representatives in Fact
 \ref{cosetminus} and the proof goes similar to above with
 appropriate changes.

\end{proof}

Now we illustrate the above proof with the example.

\bex

Consider the set $T[(2,1),1,+]$ and show explicitly that
$\rho^{(2,1),1,+}\cong
 S^{(2,1)}\uparrow_{S_3^1}^{D_3}$. The minimal length coset
 representatives are
 \begin{align*}
 &\sigma_1=e=(1\,\,2\,\,3\,)\,\,, \\
&\sigma_2=s_0=(-2\,\,-1\,\,3\,)\,\,,\\
&\sigma_3=s_2s_0=(-3\,\,-1\,\,2\,)\,\,,\\
&\sigma_4=s_1s_2s_0=(-3\,\,-2\,\,1\,)\,\,.
\end{align*}
The set $SYT(\lambda)=\{T_1,T_2\}$ is
$$
T_1=\young(12,3)\quad\textrm{and}\quad
T_2=T_1^{s_2}=\young(13,2)\,\,\,.
$$
We disregard the part with negative contents and the content of
the lower left box in $T_1$ and $T_2$ is 1.\\
 Then
$$
T_1^{s_0}=\young(\mtwo\mone,3)\quad\quad\quad
T_2^{s_0}=\young(\mtwo3,\mone)\,\,\,,
$$
$$
T_1^{s_2s_0}=\young(\mthree\mone,2)\quad\quad\quad
T_2^{s_2s_0}=\young(\mthree2,\mone)\,\,\,,
$$
$$
T_1^{s_1s_2s_0}=\young(\mthree\mtwo,1)\quad\quad\quad
T_2^{s_1s_2s_0}=\young(\mthree1,\mtwo)\,\,\,.
$$
So
$$
\rho_{\sigma_2}(T_1)={1\over 5}T_1+{\sqrt{24}\over
5}T_1^{s_0}\qquad,\qquad\rho_{\sigma_2}(T_2)={1\over
3}T_2+{\sqrt{8}\over 3}T_2^{s_0}
$$

$$
\rho_{\sigma_3}(T_1)={1\over
{10}}T_1+{\sqrt{3}\over{10}}T_1^{s_2}+{\sqrt{24}\over
{15}}T_1^{s_0}+{\sqrt{192}\over {15}}T_1^{s_2s_0}\,\,,
$$
$$
\rho_{\sigma_3}(T_2)=-{1\over 6}T_2+{\sqrt{24}\over
6}T_2^{s_2}-{\sqrt{8}\over {15}}T_2^{s_0}+{\sqrt{192}\over
{15}}T_2^{s_2s_0}\,\,,
$$

$$
\rho_{\sigma_4}(T_1)={1\over
{10}}T_1-{\sqrt{3}\over{10}}T_1^{s_2}+{\sqrt{24}\over
{15}}T_1^{s_0}+{{2\sqrt{3}}\over{15}}T_1^{s_2s_0}+{2\over\sqrt{5}}T_1^{s_1s_2s_0}\,\,,
$$

$$
\rho_{\sigma_4}(T_2)={1\over 6}T_2+{\sqrt{24}\over
6}T_2^{s_2}+{\sqrt{8}\over
{15}}T_2^{s_0}+{{2\sqrt{3}}\over{15}}T_2^{s_2s_0}+{2\over\sqrt{5}}T_2^{s_1s_2s_0}\,\,.
$$

Here one can easily see that each $\rho_{\sigma_i}(T_j)$ is
linearly independent of all other $\rho_{\sigma_{i'}}(T_{j'})$ and
therefore
$$
\bigoplus_{i=1}^{4}Span\{\rho_{\sigma_i}(T_1),\rho_{\sigma_i}(T_2)\}=Span\{T[(2,1),1,+]\}
$$
which together with the fact that the $\rho$-action of $D_3$ by
the construction permutes the subspaces
$Span\{\rho_{\sigma_i}(T_1),\rho_{\sigma_i}(T_2)\}$ , $1\leq i\leq
4$ , implies that $\rho=\rho^{(2,1),1,+}\cong
 S^{(2,1)}\uparrow_{S_3^1}^{D_3}$.

\eexa

\subsection{Decomposition of induced representation into minimal AY
representations}\label{decind}
 In this section we prove the main
result

\subsubsection{Decomposition Rule}

\begin{thm}\label{decomp}

Let $\lambda$ be a straight or skew shape with $n$ boxes and
$m\in\mathbb N$. Then
\begin{itemize}
\item[(1)]  $\rho^{\lambda,\boxtimes}\oplus\rho^{\lambda,\cdot
|\cdot,+}\oplus\rho^{\lambda,\div,+}\cong\rho^{\lambda,m,+}$

\item[(2)]  $\rho^{\lambda,\boxtimes}\oplus\rho^{\lambda,\cdot
|\cdot,-}\oplus\rho^{\lambda,\div,-}\cong\rho^{\lambda,m,-}$

\end{itemize}

\end{thm}

Note that, when the set $T[\lambda,\cdot |\cdot,\pm]$ or
$T[\lambda,\div,\pm]$ is empty, then the representation
$\rho^{\lambda,\cdot |\cdot,\pm}$ or $\rho^{\lambda,\div,\pm}$,
respectively, is the zero module.

  Notice also that in Theorem~\ref{decomp} it is enough to deal with $m=1$ because it follows
from Theorem~\ref{ind} that representations $\rho^{\lambda,m_1,+}$
and $\rho^{\lambda,m_2,+}$ for any $m_1,m_2\in\mathbb N$ are
isomorphic.
\subsubsection{Proof of Theorem \ref{decomp}}

First we will prove the following lemma:

\blem\label{dlem} There exist matrix functions
$$
g_i\,:\,\mathbb C\to\mathcal M_d(\mathbb C)
$$
for $i=0,1,2,...,n-1$ and $d=2^{n-1}f^\lambda$ such that
$$
g_i(m)=\rho^{\lambda,m,+}_{s_i}\qquad\qquad m=0,1,2,...
$$
The entries of $g_i(x)$ are analytic single valued functions on
$\{x\in\mathbb C\,:\,Re\,x>-1\}$ and are continuous from the right
at $x=-1$ as functions of real variable. The same holds for
$\rho^{\lambda,m,-}$.\elem

 The following obvious observation is crucial for the proof of
Lemma~\ref{dlem}:

\bobs\label{contperm} If $(c_1,c_2,...,c_n)$ is the content vector
of a D-Young tableau $T$, then the content vector of $T^{\pi}$ is
$\left(c_{\pi^{-1}(1)},c_{\pi^{-1}(2)},...,c_{\pi^{-1}(n)}\right)$.
\eobs

\noindent {\bf Proof of Lemma~\ref{dlem}.} Take some standard
D-Young tableau $Q=\lceil T\rfloor_m$ (for the definition of
$\lceil T\rfloor_m$ see Lemmas~\ref{T0} and~\ref{Tm}) and consider
its half of shape $\lambda$ with nonnegative contents denoted by
$T$. If for a certain $\pi\in D_n$, $T^{\pi}$ is standard, then
$Q^{\pi}=\left(\lceil T\rfloor_m\right)^\pi=\lceil T^\pi\rfloor_m$
is also standard for $m=1,2,3,...$; however, for $m=0$, we can
have a situation when $T^{\pi}$ is standard while
$Q^{\pi}=\left(\lceil T\rfloor_0\right)^\pi$ is not standard. In
such a case, $T^{\pi}$ is a half of shape $\lambda$ with
nonnegative contents of the standard D-Young tableau $\lceil
T^\pi\rfloor_0$ from some other subset of $T[\lambda,0,\pm]$, not
from that which contains $Q$. For example,
$$
Q=\young(:\mthree12,\mtwo\mone3)\in T[(2,1),\cdot |\cdot,+]\qquad
T=\young(12,3)
$$
and take $\pi=\begin{pmatrix}
-3 &-2 &-1 &1 &2 &3 \\
-2 &1   &3   &-3 &-1 &2
\end{pmatrix}\in D_3$. Then
$ T^{\pi}=\young(\mthree\mone,2) $ is standard while
$Q^{\pi}=\left(\lceil
T\rfloor_0\right)^\pi=\young(:\mtwo\mthree\mone,132)\notin
T[(2,1),\cdot |\cdot,+]$ because is not standard. Here $ T^{\pi}$
is the half of shape $\lambda=(2,1)$ of nonnegative contents for
the standard D-Young tableau $\lceil
T^\pi\rfloor_0=\young(:\mthree\mone,:\pmtwo,13)\in
T[(2,1),\boxtimes]$.

If we want to prove Theorem~\ref{decomp}(1) we take a standard
Young tableau $T$ with $n$ boxes of shape $\lambda$ filled by
numbers $1,2,3,...,n$ (the number of negative entries is zero,
thus even). Let
$$
c=cont\left(\lceil T\rfloor_0\right)=(c_1,c_2,...,c_n)
$$
and define for $x\in\{-1\}\bigcup\,\{x\in\mathbb C\,:\,Re\,x>-1\}$
the vector function
$$
[c](x)=\left(\left([c](x)\right)_1,\left([c](x)\right)_2,...,\left([c](x)\right)_n\right)
=(c_1+x+1,c_2+x+1,...,c_n+x+1)\,.
$$
For $\pi\in D_n$ define also
$$
[c]^\pi(x)=\left(\left([c](x)\right)_{\pi^{-1}(1)},\left([c](x)\right)_{\pi^{-1}(2)},...,
\left([c](x)\right)_{\pi^{-1}(n)}\right)\,.
$$
Note that $[c](x)=[c]^e(x)$ where $e\in D_n$ is the identity
element. Also note that when $\pi^{-1}(j)=-i$ for some $1\leqslant
i,j\leqslant n$, then we have
$$
\left([c](x)\right)_{\pi^{-1}(j)}=\left([c](x)\right)_{-i}=-\left([c](x)\right)_i=-(c_i+x+1)\,.
$$
Clearly, by the construction of $\lceil T\rfloor_m$ presented in
proofs of Lemmas~\ref{T0} and~\ref{Tm}, each entry of $c(\lceil
T\rfloor_m)$ is greater by 1 than the corresponding entry of
$c(\lceil T\rfloor_{m-1})$ i.e.
$$
c(\lceil T\rfloor_m)=c(\lceil T\rfloor_{m-1})+(1,1,...,1)\,.
$$
Note that this holds for $m=1$ too.

 These functions $[c](x)$ are
defined in such a way that substitutions $x=-1,0,1,2,...$ in
$[c](x)$ lead to content vectors of our tableaux $\lceil
T\rfloor_m$ with $m\in\{0\}\bigcup\mathbb N$, i.e.,
$$
[c](m-1)=c(\lceil T\rfloor_m)\,;
$$
and therefore by definition of $[c]^\pi(x)$ and
Observation~\ref{contperm}, we have
$$
[c]^{\pi}(m-1)=c(\lceil T^{\pi}\rfloor_m)\,\,.
$$
Now define
$$
[f](x)=\der[c](x)\,\,\,,\,\,\,[f]^\pi(x)=\der[c]^\pi(x)\,.
$$
The following is a most {\bf important observation} of this proof:
for integer $x=m-1$ where $m=0,1,2,...$, we have
$$
[f](m-1)=\der[c](m-1)=\der c(\lceil T\rfloor_m)=f(\lceil
T\rfloor_m)\,\,,
$$
$$
[f]^\pi(m-1)=\der[c]^{\pi}(m-1)=\der c(\lceil
T^{\pi}\rfloor_m)=f(\lceil T^\pi\rfloor_m)\,\,.
$$
Now consider the space spanned by the set $R$ of standard tableaux
of shape $\lambda$ with even number of negative entries
$$R=\{P\,:\,P\,\,\textrm{is standard and}\,\,P=T^\pi\,\,\textrm{ for some}\,\,\pi\in
D_n\}$$ and define for each $i=0,1,2,...,n-1$ a linear operator on
$Span(R)$ by its action on the basis set $R$:
$$
g_i(x+1)(T^\pi)=\frac{1}{\left([f]^{\pi}(x)\right)_i}\,T^\pi+
\sqrt{1-\frac{1}{\left([f]^{\pi}(x)\right)_i^2}}\,T^{s_i\pi}\,\,\,,
$$
where $\left([f]^{\pi}(x)\right)_i$ denotes the $i$-th entry of
the vector $[f]^{\pi}(x)$ and the branch of square root is
$\sqrt1=1$. Observe that when $T^{s_i\pi}$ is not standard, then
$\left([f]^{\pi}(x)\right)_i=\pm1$, so we shall never exit
$Span(R)$. The matrices of these operators with respect to the
basis $R$ are matrix functions $g_i(x+1)$. By the construction,
these are exactly the representation matrices of the generators
$\rho^{\lambda,m,+}_{s_i}$ for $x=m-1$. This holds since the
generator matrices of the representation $\rho^{\lambda,m,+}$ are
completely determined by entries of derived content vectors of
tableaux of the set $T[\lambda,m,+]=\{\lceil
T^{\pi}\rfloor_m\,:\,T^{\pi}\,\,\textrm{is standard}\}$.

In order to complete the proof of Lemma~\ref{dlem} we have to
prove the following

 \bcl\label{roweldef} For $i=0,1,2,...,n-1$, mappings $
g_i(x+1)$ are well-defined for $x\in\mathbb C$ with $Re\,x>-1$.
The coefficients $\frac{1}{\left([f]^{\pi}(x)\right)_i}$ and
$\sqrt{1-\frac{1}{\left([f]^{\pi}(x)\right)_i^2}}$ are single
valued and analytic. \ecl
\begin{proof}
We must show that $|\left([f]^{\pi}(x)\right)_i|>1$ for
$x\in\mathbb C$ with $Re\,x>-1$. By definition,
$$
\left([f]^{\pi}(x)\right)_0=\left([c]^{\pi}(x)\right)_1+\left([c]^{\pi}(x)\right)_2=
\left([c](x)\right)_{\pi^{-1}(1)}+\left([c](x)\right)_{\pi^{-1}(2)}\,,
$$
$$
\left([f]^{\pi}(x)\right)_i=\left([c]^{\pi}(x)\right)_{i+1}-\left([c]^{\pi}(x)\right)_i=
\left([c](x)\right)_{\pi^{-1}(i+1)}-\left([c](x)\right)_{\pi^{-1}(i)}\,
$$
for $i=1,2,...,n-1$. Obviously, if $\pi^{-1}(1)$ and $\pi^{-1}(2)$
have opposite signs, then $ \left([f]^{\pi}(x)\right)_0$ is
constant, i.e. does not depend on $x$, and if
$\pi^{-1}(i+1)\pi^{-1}(i)>0$, then $\left([f]^{\pi}(x)\right)_i$
does not depend on $x$. If $\pi^{-1}(1)\pi^{-1}(2)>0$ or
$\pi^{-1}(i+1)\pi^{-1}(i)<0$, then for $i=0,1,2,...n-1$
$$
|\left([f]^{\pi}(x)\right)_i|=|c_{i_1}+c_{i_2}+2x+2|\,,
$$
where $i_1=|\pi^{-1}(i+1)|$ and $i_2=|\pi^{-1}(i)|$ for
$i\geqslant 1$, or $i_1=|\pi^{-1}(1)|$ and $i_2=|\pi^{-1}(2)|$ for
$i=0$. So,
$$
|\left([f]^{\pi}(x)\right)_i|=\sqrt{(c_{i_1}+c_{i_2}+2Re\,x+2)^2+(2Im\,x)^2}
\geqslant|c_{i_1}+c_{i_2}+2Re\,x+2|\,.
$$
Numbers $c_{i_1}$, $c_{i_2}$ are nonnegative because they entries
of the content vector of $\lceil T\rfloor_0\in T[\lambda,0,+]$ and
at most one of them may be equal to zero. Hence
$c_{i_1}+c_{i_2}\geqslant 1$. By the assumption of our claim,
$Re\,x>-1$, and therefore $2Re\,x+2>0$. Thus, $
|\left([f]^{\pi}(x)\right)_i| \geqslant
c_{i_1}+c_{i_2}+2Re\,x+2>1$, and we are done.
\end{proof}
\noindent The proof of Lemma~\ref{dlem} is completed.\qed

\noindent For the further arguments we need the following
Carlson's theorem:
\begin{thm}\label{CarTh} {\bf Carlson's theorem.} If $f(z)$ is regular and of
exponential type in the half plane $x\geqslant 0$ (i.e. it $f$ is
analytic in the open half plane $x>0$ and continuous in $x\geq 0$)
and $h\left({\pi\over 2}\right)+h\left(-{\pi\over 2}\right)<2\pi$,
then $f(z)\equiv 0$ if $f(n)=0$, $n\in\mathbb Z_{+}$, where
$h(\theta)$ is the indicator function of $f$:
$$
h(\theta)=\limsup_{r\to+\infty}\frac{\log|f(re^{i\theta})|}{r}
$$
\end{thm}
 (See \cite{CT} for the details about the indicator function
$h(\theta)$ and the proof of the Carlson's theorem.)

\bigskip
\noindent {\bf Proof of Theorem~\ref{decomp}.} Now we give two
rather similar arguments such that each of them completes the
proof of Theorem~\ref{decomp}.

\medskip \noindent {\bf Argument (1)}  Take some $w=s_{i_t}\cdots
s_{i_2}s_{i_1}\in D_n$ and denote $$\chi_w(x)=Trace
(g_{i_t}(x+1)\cdots g_{i_2}(x+1)g_{i_1}(x+1))\,.$$ If we fix $w$,
then $\chi_w$ is a function of $x$ and it is a polynomial in
several expressions of the form
$$
\pm\frac{1}{c_{i}+c_{j}+2x+2}\qquad,\qquad\sqrt{1-\frac{1}{(c_{i}+c_{j}+2x+2)^2}}\,,
$$
where $c_{i}$, $c_{j}$ are some entries of the content vector of
$\lceil T\rfloor_0$, as we have seen in the proof of
Claim~\ref{roweldef}. Obviously, there exist real positive
constants $R$ and $d$ such that
$$
\forall x=re^{i\theta} \in\mathbb C\,\,:\,\,|x|=r\geqslant
R\,\,\Longrightarrow |\chi_w(x)|=|\chi_w(re^{i\theta})|\leqslant
r^d\, .
$$
 Notice that $\chi_w(m-1)$ (with
$m=1,2,3,...$) is the value of character of the representation
$\rho^{\lambda,m,+}$. We already know by Theorem~\ref{ind} that
all $\rho^{\lambda,m,+}$ for $m=1,2,3,...$ are isomorphic and
therefore $\chi_w(m-1)=const$ for $m=1,2,3,...$. Considering
$\chi_w(x)$ as a function of $x$ and choosing for square roots the
branch with $\sqrt1=1$ (the branching points of our square roots
are at the left of zero (indeed at the left of $-1$) as we showed
in Claim~\ref{roweldef}, and so we can take the branches which are
analytic in the right half plane), we see that $\chi_w(x)$
satisfies conditions of Carlson's theorem (see
Theorem~\ref{CarTh}):
\begin{align*}
\forall\theta
\,:\,h(\theta)&=\limsup_{r\to+\infty}\frac{\log|\chi_w(re^{i\theta})|}{r}\leqslant
\limsup_{r\to+\infty}\frac{\log
\left(r^d\right)}{r}=0\,\Longrightarrow\,\\
&\Longrightarrow\,h\left({\pi\over 2}\right)+h\left(-{\pi\over
2}\right)=0<2\pi\,\,.\nonumber
\end{align*}

Hence, by Carlson's theorem, $\chi_w(x)$ must be constant for any
number with positive real part, in particular for real $x>0$.
Moreover, by the uniqueness of analytic extension, it must be
constant also for the numbers with real part between $-1$ and 0.
So $\chi_w(x)=const$ for $\{x\in\mathbb C\,:\,Re\,x>-1\}$ Because
at the point $x=-1$ the function $\chi_w(x)$ is continuous from
the right function of real variable, we have
$\chi_w(-1)=\chi_w(x)=const$. But at $x=-1$, the value
$\chi_w(-1)$ is the value of character of the representation
$\rho^{\lambda,0,+}$ calculated at $w$. The element $w$ is an
arbitrary element of our group $D_n$, so the character of
$\rho^{\lambda,0,+}$ is the same as the character of
$\rho^{\lambda,m,+}$ which implies that
$$
\rho^{\lambda,0,+}=\rho^{\lambda,\boxtimes}\oplus\rho^{\lambda,\div,+}\oplus\rho^{\lambda,\cdot
|\cdot,+} \cong\rho^{\lambda,m,+}\,\,,
$$ and the first statement of Theorem~\ref{decomp} is proved.

\bigskip
Here is another argument which also proves the first statement of
Theorem~\ref{decomp}.\\
\noindent {\bf Argument (2)}  The group $D_n$ has $(n^2+n)/2$
Coxeter relations:
$$
(s_i s_j)^{m_{ij}}=1\,\,,\,\,0\leqslant i\leqslant j\leqslant
n-1\,,
$$
where $m_{ii}=1$ and $m_{ij}=2$ or $m_{ij}=3$ for $i<j$.

According to these relations, we introduce $(n^2+n)/2$ matrix
functions of $x$ for $0\leqslant i\leqslant j\leqslant n-1$:
$$
A_{ij}(x)=(g_i(x+1)g_j(x+1))^{m_{ij}}\,.
$$
For $x=-1,0,1,2,3,...$, the matrices $g_i(x+1)$ (for
$i=0,1,2,...,n-1$) are generator matrices of representations
$\rho^{\lambda,0,+}$, $\rho^{\lambda,1,+}$,
$\rho^{\lambda,2,+}$,..., and therefore the generator matrices
$\rho_{s_i}^{\lambda,0,+}$, $\rho_{s_i}^{\lambda,1,+}$,
$\rho_{s_i}^{\lambda,2,+}$, ... must satisfy the defining
relations of the group which implies
$$
A_{ij}(x)=I\qquad\textrm{for}\,\,\,0\leqslant i\leqslant
j\leqslant n-1\,\,\,\textrm{and}\,\,\,x=-1,0,1,2,3,...
$$
where $I$ is the identity matrix. But the entries of $A_{ij}(x)$
are polynomials in several expressions of the form
$$
\pm\frac{1}{c_{i}+c_{j}+2x+2}\qquad,\qquad\sqrt{1-\frac{1}{(c_{i}+c_{j}+2x+2)^2}}\,,
$$
where $c_{i}$, $c_{j}$ are some entries of the content vector of
$\lceil T\rfloor_0$, as we have seen in the proof of Claim
\ref{roweldef}. By the same argument (which involves the Carlson's
theorem) as above (in Argument (1)), we conclude that
$$
A_{ij}(x)=I\qquad\textrm{for}\,\,\,0\leqslant i\leqslant
j\leqslant n-1\,\,\,\textrm{and for any real}\,\,\,x\geqslant
-1\,\,,
$$
not only for integer $x=-1,0,1,2,3,...$  as before. This means
that the matrices $g_i(x+1)$ are generator matrices of
representations not only for integer $x=-1,0,1,2,3,...$  but for
any real $x\geq -1$. Denote this representation $\rho^{(x+1)}$. As
above, take some fixed $w\in D_n$ and denote by $\chi_w(x)$ the
character of the representation $\rho^{(x+1)}$ evaluated at $w$.
The character is a polynomial in the entries of representation
matrices, so in our case $\chi_w(x)$ is a polynomial in several
expressions of the form
$$
\pm\frac{1}{c_{i}+c_{j}+2x+2}\qquad,\qquad\sqrt{1-\frac{1}{(c_{i}+c_{j}+2x+2)^2}}\,,
$$
and therefore $\chi_w(x)$ is a continuous function of $x\geq -1$.
By discreteness of the character values and continuity
$\chi_w(x)$, we have $\chi_w(x)=const$ which implies that all the
representations $\rho^{(x+1)}$ for $x\geqslant -1$ are isomorphic.
In particular, substituting $x=-1,0,1,2,3,...$ we get
$$
\rho^{\lambda,\boxtimes}\oplus\rho^{\lambda,\div,+}\oplus\rho^{\lambda,\cdot
|\cdot,+} \cong\rho^{\lambda,m,+}\,\,,
$$ and the first statement of Theorem~\ref{decomp} is proved.

 In order to prove the second statement of Theorem~\ref{decomp}, we
have to start with standard tableau $T$ filled by numbers
$-1,2,3,...,n$. Everything is similar to the above with very
slight changes. The proof of Theorem~\ref{decomp} is
completed.\qed

\begin{rem}\label{clvwindrep} In the argument {\bf(2)} we did not use Theorem
\ref{ind}. Moreover, in argument~{\bf(2)} we proved that all the
representations $\rho^{(x)}$ are isomorphic and denoting
$$
\rho^{(\infty)}_{s_i}=\lim_{x\to +\infty}\rho^{(x)}_{s_i}\,\,\,\,
\,\,\,\,\textrm{for}\,\,\,i=0,1,2,...,n-1
$$
 we obtain the representation $\rho^{(\infty)}$ which is exactly the
 classical form of induced representation as given, for example,
 in~\cite{U-I}, proof of Theorem~9.3, page~32. This observation itself
can be used to give another proof of Theorem~\ref{ind} because, as
we have seen above in argument~{\bf(2)}, $\chi_w(x)=const$ and has
the same value when $x\to\infty$ which implies that
representations $\rho^{(x)}$ (and in particular
$\rho^{\lambda,m,+}$ for $m\in\mathbb N$) are isomorphic to the
representation $\rho^{(\infty)}$ obtained when $x$ tends to
infinity.
\end{rem}

\subsubsection{An Example}
Now we give an example to illustrate the above proof of
Theorem~\ref{decomp}.
\bcl\label{impex} For $m\in\mathbb N$
$$
\rho^{(3),\boxtimes}\oplus\rho^{(3),\div,+}\rho^{(3),\cdot
|\cdot,+}\cong\rho^{(3),m,+}
$$
\ecl

 In this case the set $T[(3),\cdot|\cdot,+]$ is empty. Hence,
Claim~\ref{impex} is equivalent to
$$
\rho^{(3),\boxtimes}\oplus\rho^{(3),\div,+}\cong\rho^{(3),m,+}\,\,\,.
$$
We can ``formulate" this claim graphically:
$$
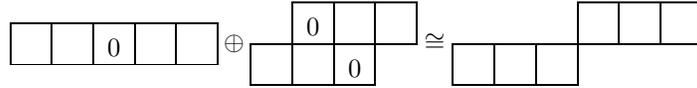
\young(\,\,0\,\,)\oplus\young(:0\,\,,\,\,0)\cong\young(:::\,\,\,,\,\,\,)
$$
Here zero shows the cells with zero content and each diagram
describes the corresponding representation:
$$
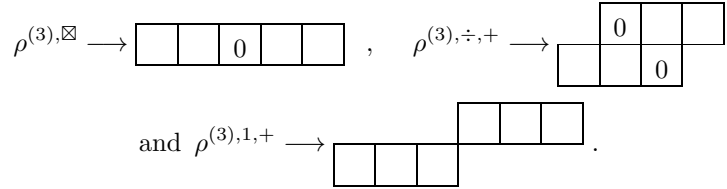
\rho^{(3),\boxtimes}\longrightarrow\young(\,\,0\,\,)\,\,\,\,\,,
\,\,\,\,\,\,\,\,\rho^{(3),\div,+}\longrightarrow\young(:0\,\,,\,\,0)
$$
$$
\textrm{and}\,\,\,\rho^{(3),1,+}\longrightarrow\young(:::\,\,\,,\,\,\,)\,\,.
$$

 Of course, this claim can be easily verified
by direct calculations because the group $D_3$ is small
(occasionally it is isomorphic to $S_4$) and its character table
is well known, but we shall follow the above general proof.

Consider four following tableaux
$$
T_1=T=\,\young(123)\,\,\,,\,\,T_2=T^{s_0}=\,\young(\mtwo\mone3)\,\,\,,
$$
$$
T_3=T^{s_2s_0}=\,\young(\mthree\mone2)\,\,\,,\,\,T_4=T^{s_1s_2s_0}=\,\young(\mthree\mtwo1)\,\,.
$$
and denote
$$
R=\{T_1,T_2,T_3,T_4\}=\left\{T^\pi\,:\,\pi\in\{e,s_0,s_2s_0,s_1s_2s_0\}\right\}\,.\eqno(*)
$$
 The half with nonnegative contents of shape $\lambda=(3)$
of standard D-Young tableau from the set $T[(3),m,+]$ for
$m\in\{0\}\cup\mathbb N$ must be one of the four tableaux from the
set $R$. \\
For each $m\in\{0\}\cup\mathbb N$ the set
$T[(3),m,+]=\left\{\lceil T_1\rfloor_m,\lceil T_2\rfloor_m,\lceil
T_3\rfloor_m,\lceil T_4\rfloor_m\right\}$. We present it in
complete details for $m=0,1$. We start with $m=0$ and we have
$T[(3),0,+]=T[(3),\boxtimes]\bigcup T[(3),\div,+]$.

 Consider $T[(3),\boxtimes]$ which consists of only one
tableau $\lceil T_1\rfloor_0$ and write down its content vector
$c(\lceil T_1\rfloor_0)$ and its derived content vector $f(\lceil
T_1\rfloor_0)$:
\begin{multline}
\lceil T_1\rfloor_0=\young(\mthree\mtwo\pmone23)\quad c(\lceil
T_1\rfloor_0)=(0,1,2)\\
\quad f(\lceil T_1\rfloor_0)=\der c(\lceil
T_1\rfloor_0)=(c_1+c_2,c_2-c_1,c_3-c_2)=\\=(0+1,1-0,2-1)=(1,1,1)\,\,.\nonumber
\end{multline}

Obviously the corresponding representation $\rho^{(3),\boxtimes}$
is the trivial representation of $D_3$ because
$\rho^{(3),\boxtimes}_{s_0}=\rho^{(3),\boxtimes}_{s_1}=\rho^{(3),\boxtimes}_{s_2}=1$.

Now consider the set $T[(3),\div,+]=\left\{\lceil
T_2\rfloor_0,\lceil T_3\rfloor_0,\lceil T_4\rfloor_0\right\}$
which is:
\begin{multline}
\lceil T_2\rfloor_0=\young(:\mtwo\mone3,\mthree12)\quad c(\lceil
T_2\rfloor_0)=(-1,0,2)\\
\quad f(\lceil T_2\rfloor_0)=\der c(\lceil
T_2\rfloor_0)=(c_1+c_2,c_2-c_1,c_3-c_2)=\\
=(-1+0,0-(-1),2-0)=(-1,1,2)\,\,,\nonumber
\end{multline}
$$
\lceil T_3\rfloor_0=\young(:\mthree\mone2,\mtwo13)\qquad f(\lceil
T_3\rfloor_0)=(1,3,-2)\,\,,
$$
$$
\lceil T_4\rfloor_0=\young(:\mthree\mtwo1,\mone23)\qquad f(\lceil
T_4\rfloor_0)=(1,-3,1)\,\,.
$$

Now we denote for brevity
$\rho^{(0)}=\rho^{(3),0,+}=\rho^{(3),\boxtimes}\oplus\rho^{(3),\div,+}$
and write down the representation matrices of generators for the
representation $\rho^{(0)}$. The upper left $1\times 1$ block in
each generator matrix corresponds to the 1-dimensional trivial
representation $\rho^{(3),\boxtimes}$ and the lower right $3\times
3$ block corresponds to the 3-dimensional representation
$\rho^{(3),\div,+}$:
$$
\rho^{(0)}_{s_0}=\begin{pmatrix}
1 &0 &0 &0\\
0 &-1 &0 &0\\
0 &0 &1 &0\\
0 &0 &0 &1
\end{pmatrix}\,\,\,,\,\,
\rho^{(0)}_{s_1}=\begin{pmatrix}
1 &0 &0 &0\\
0 &1 &0 &0\\
0 &0 &{1\over 3} &{\sqrt{8}\over 3}\\
0 &0 &{\sqrt{8}\over 3} &-{1\over 3}
\end{pmatrix}\,,\,
$$
$$
\rho^{(0)}_{s_2}=\begin{pmatrix}
1 &0 &0 &0\\
0 &{1\over 2} &{\sqrt{3}\over 2} &0\\
0 &{\sqrt{3}\over 2} &-{1\over 2} &0\\
0 &0 &0 &1
\end{pmatrix}\,.
$$

Now consider the minimal AY representation
$\rho^{(1)}=\rho^{(3),1,+}$ which is induced from the trivial
representation of $S_3^1$ to $D_3$ as it follows from
Theorem~\ref{ind}. The set $T[(3),1,+]$ is:
$$
\lceil T_1\rfloor_1=\young(:::123,\mthree\mtwo\mone)\quad c(\lceil
T_1\rfloor_1)=(1,2,3)\quad f(\lceil T_1\rfloor_1)=(3,1,1)\,\,,
$$
$$
\lceil T_2\rfloor_1=\young(:::\mtwo\mone3,\mthree12)\quad
 f(\lceil T_2\rfloor_1)=(-3,1,4)\,\,,
$$
$$
\lceil T_3\rfloor_1=\young(:::\mthree\mone2,\mtwo13)\qquad
f(\lceil T_3\rfloor_1)=(1,5,-4)\,\,,
$$
$$
\lceil T_4\rfloor_1=\young(:::\mthree\mtwo1,\mone23)\qquad
f(\lceil T_4\rfloor_1)=(1,-5,1)\,\,.
$$
The generator matrices in this case are
$$
\rho^{(1)}_{s_0}=\begin{pmatrix}
{1\over 3} &{\sqrt{8}\over 3} &0 &0\\
{\sqrt{8}\over 3} &-{1\over 3} &0 &0\\
0 &0 &1 &0\\
0 &0 &0 &1
\end{pmatrix}\,\,\,,\,\,
\rho^{(1)}_{s_1}=\begin{pmatrix}
1 &0 &0 &0\\
0 &1 &0 &0\\
0 &0 &{1\over 5} &{\sqrt{24}\over 5}\\
0 &0 &{\sqrt{24}\over 5} &-{1\over 5}
\end{pmatrix}\,\,,
$$
$$
 \rho^{(1)}_{s_2}=\begin{pmatrix}
1 &0 &0 &0\\
0 &{1\over 4} &{\sqrt{15}\over 4} &0\\
0 &{\sqrt{15}\over 4} &-{1\over 4} &0\\
0 &0 &0 &1
\end{pmatrix}\,.
$$

Define matrix functions:
$$
g_0(x+1)=\begin{pmatrix}
{1\over {2x+3}} &{\sqrt{(2x+3)^2-1}\over {2x+3}} &0 &0\\
{\sqrt{(2x+3)^2-1}\over {2x+3}} &-{1\over {2x+3}} &0 &0\\
0 &0 &1 &0\\
0 &0 &0 &1
\end{pmatrix}\,\,\,,\,
$$
$$
g_1(x+1)=\begin{pmatrix}
1 &0 &0 &0\\
0 &1 &0 &0\\
0 &0 &{1\over {2x+5}} &{\sqrt{(2x+5)^2-1}\over {2x+5}}\\
0 &0 &{\sqrt{(2x+5)^2-1}\over {2x+5}} &-{1\over {2x+5}}
\end{pmatrix}\,\,\,,
$$
$$
g_2(x+1)=\begin{pmatrix}
1 &0 &0 &0\\
0 &{1\over {2x+4}} &{\sqrt{(2x+4)^2-1}\over {2x+4}} &0\\
0 & {\sqrt{(2x+4)^2-1}\over {2x+4}} &-{1\over {2x+4}}&0\\
0 &0 &0 &1
\end{pmatrix}\,.
$$
Easy to see that substituting $x=-1,0$ we get matrices
$\rho^{(0)}_{s_i},\rho^{(1)}_{s_i}$ (for $i=0,1,2$), respectively.
The entries and the traces of the products of these matrix
functions are polynomials in $\frac{1}{2x+3}$,
$\frac{\sqrt{(2x+3)^2-1}}{2x+3}$, $\frac{1}{2x+5}$,
$\frac{\sqrt{(2x+5)^2-1}}{2x+5}$, $\frac{1}{2x+4}$,
$\frac{\sqrt{(2x+4)^2-1}}{2x+4}$. Therefore, the Carlson's theorem
is applicable, as we explained in argument~(1) of the proof of
Theorem~\ref{decomp} and both of arguments~(1) and~(2) work here
properly to show that representations $\rho^{(0)}$ and
$\rho^{(1)}$ are isomorphic.

The matrices of generators in the representation $\rho^{(\infty)}$
which was introduced in Remark \ref{clvwindrep} are
$$
\rho^{(\infty)}_{s_0}=\begin{pmatrix}
0 &1 &0 &0\\
1 &0 &0 &0\\
0 &0 &1 &0\\
0 &0 &0 &1
\end{pmatrix},\,
\rho^{(\infty)}_{s_1}=\begin{pmatrix}
1 &0 &0 &0\\
0 &1 &0 &0\\
0 &0 &0 &1\\
0 &0 &1 &0
\end{pmatrix},\,
\rho^{(\infty)}_{s_2}=\begin{pmatrix}
1 &0 &0 &0\\
0 &0 &1 &0\\
0 &1 &0 &0\\
0 &0 &0 &1
\end{pmatrix}\,.
$$
The representation $\rho^{(\infty)}$ is induced into $D_3$ from
the trivial representation of $S_3^1=\langle s_1,s_2\rangle$ and
is the classical form of induced representation as given, for
example, in~\cite{U-I}, proof of Theorem~9.3, page~32. We see that
in our case the representation $\rho^{(\infty)}$ is a permutation
representation, moreover its generator matrices are exactly as in
the natural representation of $S_4$. It is not a surprise because
$D_3\cong S_4$ and the representation which is induced from the
trivial representation is a permutation representation obtained
from the action of $D_3$ on the set $D_3/S_3^1$ which consists of
four cosets. It is a well known fact that the natural
representation of the symmetric group $S_n$ is isomorphic to the
direct sum the trivial representation $S^{(n)}$ and the
irreducible representation of $S^{(n-1,1)}$ associated with the
partition $(n-1,1)\vdash n$. In our example we have a special case
of this observation for $n=4$: the representation
$\rho^{(3),\boxtimes}$ is the trivial representation of $D_3\cong
S_4$ while the representation $\rho^{(3),\div,+}$ is isomorphic to
the irreducible representation of $S_4\cong D_3$ associated with
the partition $(3,1)\vdash 4$.

In the previous example the decomposition had only two factors (we
had there $\lambda=(3)$, the set $T[\lambda,\cdot
|\cdot,+]=\varnothing$ and the representation $\rho^{\lambda,\cdot
|\cdot,+}$ was the zero representation). I f we take
$\lambda=(2,2)$ then all three factors $\rho^{\lambda,\cdot
|\cdot,\pm}$, $\rho^{\lambda,\div,\pm}$ and
$\rho^{\lambda,\boxtimes}$ are nonzero representations. This can
be expressed graphically
$$
\young(\,0\,\,,\,\,0\,)\oplus\young(:\,\,,\,0\,,\,\,)\oplus
\young(\,\,,0\,,\,0,\,\,)\cong\young(::\,\,,::\,\,,\,\,,\,\,)
$$
Here zero shows the boxes with zero content and each diagram
describes the corresponding representation:
$$
\rho^{2^2,\cdot
|\cdot,\pm}\longrightarrow\young(\,0\,\,,\,\,0\,)\,\,\,\,\,\,\,,
\,\,\,\,\,\,\,\rho^{2^2,\boxtimes}
\longrightarrow\young(:\,\,,\,0\,,\,\,)\,\,,
$$
$$
\rho^{2^2,\div,\pm}\longrightarrow\young(\,\,,0\,,\,0,\,\,)\,\,\,\,\,\,\,,\,\,\,\,\,\,\,
\rho^{2^2,1,\pm}\longrightarrow\young(::\,\,,::\,\,,\,\,,\,\,)\,\,.
$$
Concerning this example we notice also the following: Theorem
\ref{ind} implies that the representation $\rho^{2^2,1,+}$ and
$\rho^{2^2,1,+}$ are isomorphic to the representations induced to
$D_4$ from the irreducible two dimensional representation of
$S_4^1$ and $S_4^0$ respectively, associated with the Young
diagram $2^2=(2,2)=\young(\,\,,\,\,)$ .

\subsubsection{Remarks}
We finish this section with several remarks.
\begin{rem} It is possible to avoid using Carlson's theorem (which is a
highly nontrivial result from the theory of analytic functions) in
arguments~{\bf (1)} and~{\bf (2)} of the proof of
Theorem~\ref{decomp}: if we take a stochastic normalization
instead of symmetric (see Remark~\ref{normalization}), then the
characters will be rational functions of $x$ (without square
roots) and a rational function which has infinitely many zeros is
zero because a polynomial can't have infinitely many zeros.
\end{rem}

\bigskip

 There is a natural question: when the
$D_n$-representations $\rho^{\lambda,m,+}$ and
$\rho^{\lambda,m,-}$ are isomorphic? By Theorem \ref{ind} we have
$\rho^{\lambda,m,+}\cong S^{\lambda}\uparrow_{S_n^1}^{D_n}$ and
$\rho^{\lambda,m,-}\cong S^{\lambda}\uparrow_{S_n^0}^{D_n}$. The
subgroups $S_n^1$ and $S_n^0$ are always conjugate in $B_n$. If
$n$ is odd, then $S_n^1$ and $S_n^0$ are also conjugate as
subgroups in $D_n$ (the element $\pi$ such that $\pi(1)=1$ and
$\pi(i)=-i$ for $i=2,3,...n$ belongs to $D_n$ when $n$ is odd and
$\pi^{-1}s_1\pi=s_0$, $\pi^{-1}s_i\pi=s_i$ for $i=2,3,...n-1$)
which easily implies that the induced representations are
isomorphic. When $n$ is even the representations
$\rho^{\lambda,m,+}$ and $\rho^{\lambda,m,-}$ can be non
isomorphic. The irreducible representations of $B_n$ are indexed
by ordered pairs of partitions $(\alpha,\beta)$ such that
$|\alpha|+|\beta|=n$. Restricting irreducible representations
$(\alpha,\beta)$ and $(\beta,\alpha)$ with $\alpha\neq\beta$ from
$B_n$ to $D_n$ we get the same irreducible representation of $D_n$
and therefore irreducible representations of $D_n$ are indexed by
unordered pairs of partitions $\{\alpha,\beta\}$ such that
$|\alpha|+|\beta|=n$, $\alpha\neq\beta$. For even $n$ there are
two non isomorphic irreducible representations of $D_n$
$\{\alpha,\alpha\}^{+}$ and $\{\alpha,\alpha\}^{-}$ where
$\alpha\vdash{n\over 2}$ which split from the irreducible
representation of $B_n$ associated with $(\alpha,\alpha)$. It can
be easily shown that
$$
1\uparrow_{S_n^0}^{B_n}\cong
1\uparrow_{S_n^1}^{B_n}\cong\bigoplus_{k=0}^{n}((n-k),(k))
$$
and therefore
$$
\rho^{(n),m,+}\cong 1\uparrow_{S_n^1}^{D_n}\cong
\bigoplus_{k=0}^{{n\over
2}-1}\{(n-k),(k)\}\oplus\left\{\left({n\over
2}\right),\left({n\over 2}\right)\right\}^{+}
$$
and
$$
 \rho^{(n),m,-}\cong 1\uparrow_{S_n^0}^{D_n}\cong
\bigoplus_{k=0}^{{n\over
2}-1}\{(n-k),(k)\}\oplus\left\{\left({n\over
2}\right),\left({n\over 2}\right)\right\}^{-}
$$
and therefore $\rho^{(n),m,+}\ncong\rho^{(n),m,-}\,$ .
\begin{rem}
Using the character table of $D_4$ one can see that the
representation $\rho^{(2^2),\cdot |\cdot,-}$ whose representation
space is spanned by
$$
\young(\mfour\mthree\mtwo1,\mone234)\,\,\,,\,\,\,
\young(\mfour\mthree\mone2,\mtwo134)\,\,\,,\,\,\,
\young(\mfour\mtwo\mone3,\mthree124)\,\,\,
$$
and whose generator matrices are
$$\rho_{s_0}^{(2^2),\cdot |\cdot,-}=\begin{pmatrix}
1 & 0 & 0  \\
0 & 1 & 0  \\
0 & 0 & -1 \\

\end{pmatrix}\,\,\,,\,\,\,
\rho_{s_1}^{(2^2),\cdot |\cdot,-}=\begin{pmatrix}
\frac{1}{3}  & \frac{\sqrt{8}}{3} &0 \\
\frac{\sqrt{8}}{3}  & -\frac{1}{3} &0 \\
0 &0 &1\\
\end{pmatrix}
$$
$$\rho_{s_2}^{(2^2),\cdot |\cdot,-}=\begin{pmatrix}
1 &0 &0\\
0 &-\frac{1}{2} &\frac{\sqrt{3}}{2}    \\
0 & \frac{\sqrt{3}}{2} &\frac{1}{2} \\

\end{pmatrix}\,\,\,,\,\,\,
\rho_{s_3}^{(2^2),\cdot |\cdot,-}=\begin{pmatrix}
1 & 0 & 0  \\
0 & 1 & 0  \\
0 & 0 & -1 \\

\end{pmatrix}
$$
is isomorphic to the irreducible threedimensional split
representation $\{(2),(2)\}^{-}$ while the representation
$\rho^{(2^2),\cdot |\cdot,+}$ whose representation space is
spanned by
$$\young(\mfour\mthree12,\mtwo\mone34)\,\,\,,\,\,\,
\young(\mfour\mtwo13,\mthree\mone24)\,\,\,,\,\,\,
\young(\mfour\mthree\mtwo\mone,1234)$$ is isomorphic to another
irreducible threedimensional split representation
$\{(2),(2)\}^{+}$. Of course, $\rho^{(2^2),\cdot |\cdot,+}\ncong
\rho^{(2^2),\cdot |\cdot,-}$ because $\{(2),(2)\}^{+}\ncong
\{(2),(2)\}^{-}$. Similarly, one can verify directly that
$\rho^{(2^2),\div,-}\cong\{(1^2),(1^2)\}^{-}$ and
$\rho^{(2^2),\div,+}\cong\{(1^2),(1^2)\}^{+}$.
\end{rem}
The following conjecture seems to be true:
\bcon Let $n\in \mathbb
N$ be even. Then
$$
\rho^{\left(\left({n\over 2}\right)^2\right),\cdot |\cdot,-}\cong
\left\{\left({n\over 2}\right),\left({n\over
2}\right)\right\}^{-}\,\,,\,\,\rho^{\left(\left({n\over
2}\right)^2\right),\cdot |\cdot,+}\cong \left\{\left({n\over
2}\right),\left({n\over 2}\right)\right\}^{+}\,,
$$
$$
\rho^{\left(2^{n\over 2}\right),\div,-}\cong
\left\{\left(1^{n\over 2}\right),\left(1^{n\over
2}\right)\right\}^{-}\,\,,\,\, \rho^{\left(2^{n\over
2}\right),\div,+}\cong \left\{\left(1^{n\over
2}\right),\left(1^{n\over 2}\right)\right\}^{+}\,.
$$

 \econ
\subsubsection{Why we need the negative contents?}\label{necessneg}

The representations $\rho^{\lambda,m,\pm}$ (with $m\in\mathbb N$),
$\rho^{\lambda,\cdot |\cdot,\pm}$, $\rho^{\lambda,\div,\pm}$,
$\rho^{\lambda,\boxtimes}$ are completely determined by the
derived content vector $f=\der c$ which itself is completely
determined by the half tableau of shape $\lambda$ which consists
of boxes with nonnegative contents(when the shape is placed; i.e.
$m$ is given or the conditions $\boxtimes$ or $\cdot|\cdot$,
$\div$ hold). The natural question is: why we need the half with
non positive contents? The answer is: we want the minimal AY cell
$\K^{c}$ which corresponds to the minimal AY representation
$\rho^c$ to be described bijectively by the set of standard
D-Young tableaux and to provide this we need
$$
f_i(T)=\pm 1\Longleftrightarrow T^{s_i}\,\,\,\textrm{is not
standard.}
$$
Consider, for example, $T=\young(123)$ where the box $\young(1)$
has the zero content. Then the derived content vector
$$
f(T)=(f_0,f_1,f_2)=(c_1+c_2,c_2-c_1,c_3-c_2)=(0+1,1-0,2-1)=(1,1,1)\,\,.
$$
We have $f_0=1$ and therefore the tableaux $T^{s_0}$ must be not
standard, but $T^{s_0}=\young(\mtwo\mone3)$ is standard.

However, if we add to $T$ the part with non positive contents and
let $T=\young(\mthree\mtwo\pmone23)$ we have
$T^{s_0}=\young(\mthree1\pmtwo\mone3)$ which is indeed not
standard as we wish.

When the tableaux has no boxes of zero content (i.e. in the case
of representation $\rho^{\lambda,m,\pm}$ with $m\in\mathbb N$),
then indeed the part with negative contents does not play any role
and we introduced it only in order to have a unified definition of
a D-Young tableau.

Ram in \cite{Ra} requires for standardness increasing of entries
in the same diagonal (from northwest to southeast) without
allowing two entries to occupy the same box. In the above example
it works properly:
$$
T=\young(\mthree\mtwo\mone,:::123)
$$
is standard and
$$
T^{s_0}=\young(\mthree12,:::\mtwo\mone3)
$$
has $2$ to the northwest of $-2$ in the same diagonal and hence is
not standard. However for
$$
T=\young(::::23,::\mone,:::1,\mthree\mtwo)\quad\textrm{with}\quad
f_0(T)=c_1(T)+c_2(T)=0+3=3\neq\pm1\,\,,
$$
we have
$$
T^{s_0}=\young(::::\mone3,::2,:::\mtwo,\mthree1)
$$
with $2$ to the northwest of $-2$ in the same diagonal, i.e.
$T^{s_0}$ is not ``standard" while $f_0(T)\neq\pm1$. If we allow
the numbers $\pm i$ to occupy the same box the (derived) content
vector does not change; however, we avoid here the mentioned
problem:
$$
T=\young(:::::23,:::\pmone,\mthree\mtwo)\,\,\,\,\,\textrm{and}\,\,
\,\,\,T^{s_0}=\young(:::::\mone3,:::\pmtwo,\mthree1)
$$
are both standard.

\section{Analogous constructions for $B_n$}\label{bn}

In this section we show how to obtain representations induced from
$S_n$ to $B_n$ by the construction which is similar to presented
above construction for $D_n$.

The tableaux in this case will be with semi integer contents:
\bde\label{semiintcont} For a (straight or a skew) tableau $T$
denote
$$c_k:=j-i+\frac{1}{2}$$
 where $k$ is the entry in row $i$ and
column $j$ of $T$. Call $cont(T):=(c_1,\dots,c_n)$ the {\em
content vector} of $T$, and call $\der cont(T)$ the {\em derived
content vector} of $T$. Below we sometimes shall denote  for
brevity $cont(T)$ as $c(T)$. \ede

\bde\label{byoungtab}
 Let $\lambda$ be a diagram of skew shape. Define a $B$-Young
tableau of shape $\lambda$ to be a filling of $\lambda$ by the
$2n$ numbers $\pm 1,\pm 2,\ldots,\pm n$ in such a way that
$c_{-i}=-c_i$ for $1 \leq i \leq n$ where $c_i$ are semi integer
contents defined above in Definition \ref{semiintcont}. A
$B$-Young tableau is called {\it standard} if the numbers are
increasing in rows and in columns. \ede

The additional generator $s_0=(-1,1)$ of $B_n$ differs from $s_0$
of $D_n$ (see sections \ref{bdef} and \ref{ddef}) and its simple
root is $e_1$. \bde\label{bdercontvect} For the B-Young tableau
$T$ with its content vector $c(T)$ with semi integer entries we
define the derived content vector
\begin{multline}
f(T)=(f_0(T),f_1(T),f_2(T),...,f_{n-1}(T))=\\
=\Delta c(T)=(2c_1,c_2-c_1,c_3-c_2,...,c_n-c_{n-1})\,.\nonumber
\end{multline}
\ede

\subsection{Construction and examples of standard B-Young tableaux}

We deal here only with standard B-Young tableaux which have at
most one box with content ${1\over 2}$ (and therefore, because of
the condition $c_{-i}=-c_i$, they also have at most one box with
content $-{1\over 2}$).

Each such tableau can divided into two parts one with positive
contents and another one with negative contents. Similarly to
construction for $D_n$ presented above in section \ref{dnconst} we
shall denote the shape of boxes with positive contents as
$\lambda$.

\subsubsection{Standard B-Young tableaux without boxes of
content $\pm{1\over 2}$.}

\bde Let $\lambda$ be a straight or skew shape with $n$ boxes and
let $m\in\mathbb N$. Then denote by $T[\lambda,m]$ the set of all
standard B-Young tableaux whose shape of boxes with positive
contents is $\lambda$ and the smallest positive content is
$m+{1\over 2}$. \ede

For example, the set $T[1^3,1]=\{T_1,T_2,...,T_8\}$ consists of
the following tableaux:
$$
T_1=\,\young(::1,::2,::3,\mthree,\mtwo,\mone)\,\,;\,
T_2=\,\young(::\mtwo,::\mone,::3,\mthree,1,2)\,\,; \,
T_3=\,\young(::\mthree,::\mone,::2,\mtwo,1,3)\,\,;\,T_4=\,\young(::\mthree,::\mtwo,::1,\mone,2,3)\,\,;
$$
$$
T_5=\,\young(::\mone,::2,::3,\mthree,\mtwo,1)\,\,;\,T_6=\,\young(::\mtwo,::1,::3,\mthree,\mone,2)\,\,;
\,T_7=\,\young(::\mthree,::1,::2,\mtwo,\mone,3)\,\,;\,T_8=\,\young(::\mthree,::\mtwo,::\mone,1,2,3)\,\,.
$$

\bpr  Let $\lambda$ be a straight or skew shape with $n$ boxes and
let $m\in\mathbb N$. Then
$$
\# T[\lambda,m]=2^n f^{\lambda}\,\,\,.
$$
\epr
\begin{proof} Similar to the proof of Proposition \ref{inddim}.
\end{proof}
Notice that there exists an obvious natural bijection
$$T[\lambda,m]\longleftrightarrow T[\lambda,m,+]\bigcup
T[\lambda,m,-]\,\,\,\textrm{ for }\,\,\,m\in\mathbb N\,.$$

\subsubsection{Standard B-Young tableaux with a single box on
${1\over 2}$-content diagonal}

\bde Denote by $T[\lambda,\cdot | \cdot]$ the set of standard
B-Young tableaux which have a single box on ${1\over 2}$-content
diagonal (and also a single box with content $-{1\over 2}$) and
which can be divided by the \textbf{vertical} straight line into
two parts (with $n$ boxes in each part)-- one consists of the
boxes with negative contents and another of positive contents --
and $\lambda$ is the shape (straight or skew) of boxes with
positive contents after this separation.
 \ede

\bde Denote by $T[\lambda,\div]$ the set of standard B-Young
tableaux which have a single box on ${1\over 2}$-content diagonal
and which can be divided by the \textbf{horizontal} straight line
into two parts (with $n$ boxes in each part)-- one consists of the
boxes with negative contents and another of positive contents --
and $\lambda$ is the shape (straight or skew) of boxes with
positive contents after this separation. \ede

For example, $T[(2,1),\cdot | \cdot]$ consists of following eleven
tableaux:
$$
\young(::12,:\mthree3,\mtwo\mone)\quad\young(::13,:\mtwo2,\mthree\mone)\quad\young(::\mone2,:\mthree3,\mtwo1)
\quad\young(::\mone3,:\mtwo2,\mthree1)
$$
$$
\young(::\mtwo1,:\mthree3,\mone2)\quad\young(::\mtwo3,:\mone1,\mthree2)\quad\young(::\mthree1,:\mtwo2,\mone3)
\quad\young(::\mthree2,:\mone1,\mtwo3)
$$
$$
\young(::\mtwo\mone,:\mthree3,12)\quad\young(::\mthree\mone,:\mtwo2,13)\quad\young(::\mthree\mtwo,:\mone1,23)
$$
The set $T[(2,1),\div]$ consists of following five tableaux:
$$
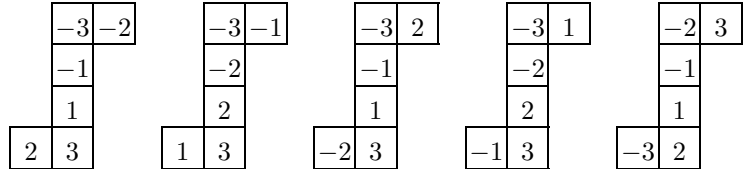
\young(:\mthree\mtwo,:\mone,:1,23)\quad\young(:\mthree\mone,:\mtwo,:2,13)\quad
\young(:\mthree2,:\mone,:1,\mtwo3)\quad\young(:\mthree1,:\mtwo,:2,\mone3)\quad
\young(:\mtwo3,:\mone,:1,\mthree2)
$$

Notice that unlike the case of $D_n$, here the sets
$T[\lambda,\cdot | \cdot]$ and $T[\lambda,\div]$ never are empty.

In the above example we see that
$$
\# T[(2,1),\cdot | \cdot]+\#T[(2,1),\div]=11+5=2^3\cdot 2=\#
T[(2,1),m]\,\,,\,\,m=1,2,...
$$
It does not surprise in view of the following proposition:

 \bpr\label{bnbij} Let $\lambda$ be a straight
or skew shape with $n$ boxes and let $m\in\mathbb N$. There exists
a natural bijection
$$
T[\lambda,\cdot | \cdot]\bigcup T[\lambda,\div]\longleftrightarrow
T[\lambda,m]
$$
\epr
\begin{proof}
Take some tableau $T\in T[\lambda,m]$ and denote by $T_{pos}$ its
half of shape $\lambda$ with positive contents and by $T_{neg}$
its half with negative contents. If the lower left box of
$T_{pos}$ is occupied by a positive entry, then we can stick
$T_{pos}$ and $T_{neg}$ by a common vertical edge to obtain a
tableau which belongs to $T[\lambda,\cdot | \cdot]$. For example,
if $\lambda=(2,1)$ and $T_{pos}=\young(\mthree\mone,2)$,
$T_{neg}=\young(:\mtwo,13)$, then we have a tableau
$\young(::\mthree\mone,:\mtwo2,13)\in T[\lambda,\cdot | \cdot]$
however sticking in the second way we obtain a tableau
$\young(:\mthree\mone,:2,:\mtwo,13)$ which is not standard.

If the lower left box of $T_{pos}$ is occupied by a negative
entry, then we can stick $T_{pos}$ and $T_{neg}$ by a common
horizontal edge to obtain a tableau which belongs to
$T[\lambda,\div]$. For example, if
$T_{pos}=\young(\mthree2,\mone)$, then we stick it to
$T_{neg}=\young(:1,\mtwo3)$ from above to obtain a standard
B-Young tableau $\young(:\mthree2,:\mone,:1,\mtwo3)\in
T[(2,1),\div]$ while sticking in the other way leads to a non
standard tableau $\young(::\mthree2,:1\mone,\mtwo3)$.

For the inverse mapping we take a tableau from $T[\lambda,\cdot |
\cdot]\bigcup T[\lambda,\div]$ and separate it by a vertical or
horizontal straight line, respectively, into two halves $T_{pos}$
and $T_{neg}$ with positive and negative contents and then put
these halves in such a way that the lower left box of $T_{pos}$
will have content $m+{1\over 2}$ (and of course the upper right
box of $T_{neg}$ will occupy the box with content
$-\left(m+{1\over 2}\right)$). This gives us a standard B-Young
tableau from the set $T[\lambda,m]$.
\end{proof}

\bde Denote $T[\lambda,0]:=T[\lambda,\cdot | \cdot]\bigcup
T[\lambda,\div]$. \ede

Reformulating Proposition \ref{bnbij}  we can say that there
exists a natural bijection between $T[\lambda,0]$ and
$T[\lambda,m]$ with $m=1,2,3,...$

\subsection{Representations arising from standard B-Young
tableaux}

As above, for a B-Young tableau $T$ and $\pi\in B_n$ we denote by
$T^{\pi}$ the B-Young tableau where each entry $i$ is replaced by
$\pi(i)$.

Take a set of standard B-Young tableaux $T[\lambda,m]$ or
$T[\lambda,\cdot |\cdot]$ or $T[\lambda,\div]$ and define a
representation $\rho^{\lambda,m}$ or $\rho^{\lambda,\cdot |\cdot}$
or $\rho^{\lambda,\div}$ by Young orthogonal form
$$\rho_{s_i}(T)=\frac{1}{f_i(T)}\,\,T+\sqrt{1-\frac{1}{f_i(T)^2}}\,\,T^{s_i}
\,\,\,\,\,\,\,\textrm{for}\,\,i=0,1,2,...,n-1\,.
$$
Here $s_0$ is $s_0$ of $B_n$ and the numbers $f_i(T)$ are the
entries of the derived content vector of B-Young tableaux $T$
introduced in Definition \ref{bdercontvect}.\\
 (It is a routine
but not so difficult work to verify that this is indeed a
representation i.e. that matrices $\rho_{s_i}$ satisfy the group
relations of $B_n$. We already know that for $i=1,2,...,n-1$ the
relations are satisfied from the classical Young orthogonal form
for $S_n$. Easy to see that
\begin{multline}
\rho_{s_0^2}(T)=\rho_{s_0}\left(\rho_{s_0}(T)\right)=\\
=\rho_{s_0}\left(\frac{1}{2c_1}\,T+\sqrt{1-\frac{1}{\left(2c_1\right)^2}}\,\,T^{s_0}\right)=\\
=\frac{1}{2c_1}\left(\frac{1}{2c_1}\,T+\sqrt{1-\frac{1}{\left(2c_1\right)^2}}\,\,T^{s_0}\right)+\\
+\sqrt{1-\frac{1}{\left(2c_1\right)^2}}\,\left(-\frac{1}{2c_1}\,T^{s_0}
+\sqrt{1-\frac{1}{\left(2c_1\right)^2}}\,\,T\right)=T
\,\,.\nonumber
\end{multline}
 Also it is rather easy to show that $\rho_{s_0}$ commutes with $\rho_{s_i}$ for $i=2,3,...,n-1$. The most
 long is the calculation which verifies that $\left(\rho_{s_0}\rho_{s_1}\right)^4=id$.)

Now we can explain why for $B_n$ we use semi integer contents. As
we mentioned above, to obtain a representation from the derived
content vector its entries must repeat the form of simple roots of
the group, therefore the first entry of the derived content vector
of B-Young tableau must be $c_1$ (or its multiple) because the
root which correspond to $s_0$ in $B_n$ is $e_1$. In the
definition of the representation $\rho$ the entry $f_0(T)=2c_1(T)$
appears in the denominator and therefore to avoid division by zero
we introduce semi integer contents $...,-{3\over 2},-{1\over
2},{1\over 2},{3\over 2},...$
 We want the representation space to be spanned by standard tableaux and so we
 need $T^{s_i}$ to be not standard if and only if
$\sqrt{1-\frac{1}{f_i(T)^2}}=0 $ which occurs if and only if
$f_i=\pm1$. This is the reason for taking $2c_1$ but not $c_1$.

The analogues of Theorems \ref{ind} and \ref{decomp} here are as
follows:
\begin{thm}\label{bnind}  Let $\lambda$ be a straight or skew shape with $n$ boxes and let
$S^{\lambda}$ denote the representation of $S_n=S_n^1=\langle
s_1,s_2,...s_{n-1}\rangle$ associated with $\lambda$ via classical
Young orthogonal form (see Theorems~\ref{t.yof}
and~\ref{t.yof-skew}) . (If $\lambda$ is a straight shape, then
$S^{\lambda}$ is an irreducible Specht module and if $\lambda$ is
a skew shape, then $S^{\lambda}$ is a skew Specht module.)
 Then for any natural $m=1,2,3,\ldots$
$$
\rho^{\lambda,m}\cong
 S^{\lambda}\uparrow_{S_n}^{B_n}
$$
\end{thm}

\begin{thm}\label{bndecomp} In the conditions of the above Theorem \ref{bnind}
$$
\rho^{\lambda,\cdot |\cdot}\oplus\rho^{\lambda,\div}
\cong\rho^{\lambda,m}
$$
\end{thm}
The proofs of Theorems~\ref{bnind} and~\ref{bndecomp} are very
similar to the proofs of Theorems~\ref{ind} and~\ref{decomp},
respectively.

A natural question is why we don't consider tableaux with semi
integer contents to obtain AY representations of $D_n$? We explain
the reason by the following example:
$$
T=\,\young(\mthree\mtwo\mone123)\qquad c(T)=\left({1\over
2}\,,\,{3\over 2}\,,\,{5\over 2}\right)
$$
For $B_3$ we have $f(T)=(2c_1,c_2-c_1,c_3-c_2)=(1,1,1)$ (this is
the trivial representation of $B_3$) and indeed for $s_0=(-1,1)$
the first entry of $f(T)$ is $f_0=2c_1=2\cdot{1\over 2}=1$ and
$T^{s_0}=\young(\mthree\mtwo1\mone23)$ is not standard. However
for $D_3$ with $s_0=(-1,2)(-2,1)$ the tableau
$T^{s_0}=\young(\mthree12\mtwo\mone3)$ is still not standard while
the entry $f_0=c_1+c_2={1\over 2}+{3\over 2}=2\neq \pm1$. We avoid
this situation because, as we already mentioned in the section
\ref{necessneg},
 we want the
minimal AY cell $\K^{f}$ which corresponds to the minimal AY
representation $\rho^f$ to be described bijectively by the set of
standard D-Young tableaux and to provide this we need
$$
f_i(T)=\pm 1\Longleftrightarrow T^{s_i}\,\,\,\textrm{is not
standard.}
$$

\noindent {\bf Acknowledgments.} I want to thank my PhD advisor
Professor Yuval Roichman for introducing me in this field and for
many helpful discussions. Also I have to express the special
gratitude to  Dr. Eli Bagno for his essential assistance in
writing sections \ref{prenot} and \ref{s.cells_Sn} of this paper.
Also I am grateful to Professor Ron M. Adin for helpful
discussions.

\end{document}